\documentclass[a4paper,10pt]{article}

\usepackage[margin=1in]{geometry}

\usepackage{amssymb}
\usepackage{amsmath}
\usepackage{amsthm}

\usepackage[T1]{fontenc}
\usepackage{newtxtext}
\usepackage[varvw]{newtxmath}

\usepackage{graphicx}
\usepackage{booktabs}
\usepackage{subcaption}
\usepackage{caption}
\usepackage{algorithm}
\usepackage{algpseudocode}
\usepackage{cite}
\usepackage{multicol}
\usepackage[bottom]{footmisc}
\usepackage{xcolor}
\usepackage{hyperref}
\newtheorem{rmk}{Remark}[section]

\title{Control of kinetic opinion dynamics in popularity-adaptive social networks}

\author{
Giacomo Albi\thanks{University of Verona, Department of Computer Science, Strada le Grazie 15, 37134 Verona (Italy). Email: \texttt{giacomo.albi@univr.it}} 
\and
Elisa Calzola\thanks{University of Ferrara, Department of Mathematics and Computer Science, Via Nicol\`o Machiavelli 30, 44121 Ferrara (Italy). Email: \texttt{elisa.calzola@unife.it}}
\and
Matteo Piu\thanks{University of Verona, Department of Computer Science, Strada le Grazie 15, 37134 Verona (Italy). Email: \texttt{matteo.piu@univr.it}}
}

\date{} 

\makeindex

\begin{document}
\newcommand{\revision}[1]{{\color{black}{#1}}}

\maketitle
\abstract{This paper presents a mathematical model for opinion dynamics in popularity-adaptive social networks, where both opinion spreading and the evolution of social media contacts depend on agents' popularity and the prominence of their views. While previous approaches accounted for the influence of popularity on opinion dynamics, we introduce a novel feedback mechanism in which opinion affects the formation of contacts. 
Within a kinetic modeling framework, we describe the evolution of the coupled dynamics of opinions and network structure, incorporating a class of control laws in order to promote interactions with popular individuals and amplify dominant opinions. Such control strategies are introduced to influence both opinion formation and connectivity, representing interventions such as awareness campaigns or moderation policies. Numerical results show how control strategies can mitigate polarization, foster consensus, or guide opinion distributions in dynamically evolving networks.}

\section{Introduction}
\label{sec:intro}
Kinetic theory has become a central mathematical framework for the description of
large systems of interacting agents, providing tractable models capable of capturing
the emergence of collective behaviors across a wide spectrum of applications,
ranging from economics and traffic dynamics to the social and biological sciences
\cite{PT13,Bellomo}. In these settings, Boltzmann-type equations allow to
derive mesoscopic descriptions starting from simple interaction rules, thus linking
microscopic decision processes to macroscopic observables. Applications of such paradigm span socio–economic models for price formation and wealth redistribution
\cite{Cor05,Par14,Cha07,T22,Burger,DPT18}, kinetic descriptions of population and
epidemic dynamics \cite{Dim21,Cha04,BDP21}, and applications in transportation engineering and crowd modelling
\cite{albi2021optimized,cristiani2014multiscale,Gun02,Pup16,ABFH19,BGQR22,Wolfram}.

Within this kinetic framework, opinion formation is a central example for interaction mechanisms, where 
repeated social exchanges give rise to dynamics well described by Boltzmann-type models, capturing consensus,
polarisation and fragmentation as studied in \cite{To06,Gal82,DT19,Dim20,fakenews}. Of particular importance is the extension of these modelling frameworks to account
for an additional layer of complexity such as individual popularity, typically quantified by the number of connections
or followers in a social networks. Hence, recent developments have focused on coupling opinion evolution with the
dynamics of the underlying social connectivity structure, see for example \cite{Ace11,Albi24,quattrociocchi2014opinion,fagioli2024opinion}.

Here we focus on kinetic model for the evolution of joint density of opinions and contacts, where each agent is represented by an
opinion $v \in [-1,1]$ and a popularity variable $c \ge 0$ measuring its social
connections. Following the approach proposed in \cite{albi2017opinion,Albi24}, we assume that dynamics of opinions depend on agents’
popularity. Furthermore, we introduce a feedback mechanism in which agents’ opinions influence the
evolution of their social connections, while the strength of opinion interactions
remains modulated by the individuals’ relative popularity. This second assumption has not been considered in previous works within this context, nevertheless it is supported by empirical evidence showing that individuals tend to form or dissolve connections based on opinion similarity. This phenomenon is widely
documented in studies on homophily, echo chambers, and opinion-driven rewiring in
online social networks, see for instance \cite{cinelli2021echo,himelboim2013tweeting,mcpherson2001birds,nguyen2020echo}.

In this modelling framework, we also investigate  the impact of control
mechanisms acting on both the contact dynamics and the opinion evolution.
Control strategies in opinion modelling have been explored at various scales see for example
\cite{during2024breaking,albi2025data,bicego2025computation}. Here we follow the
approach introduced in \cite{APZ14,albi2015kinetic}, deriving two  type of feedback controls, the first acting over the contact formation and designed
to counteract the endogenous decay of contacts by mitigating unfavourable trends in
the visibility of selected agents or groups. A second control term acts at
the opinion level, steering agents toward prescribed target positions whenever
activation criteria based on connectivity or local opinion density are satisfied.

The integration of such kinetic framework with control thus offers a rigorous tool for
designing strategies that promote consensus, limit fragmentation, or may regulate the
impact of highly connected agents in public discourse.
Overall, the finding of the present manuscript provides a modeling and computational setting for to \revision{analyze} interventions in
complex social systems where popularity, influence, and opinions co-evolve.

The paper is structured as follows. In Section \ref{sec:contactmodel}, starting from a system of indistinguishable agents, we present the microscopic model for the evolution of the contacts. Such evolution is driven both by the desire of the agents to increase their followers once they enter the network and, consequently, gain more influence, but it also depends on how popular and endorsed are their opinions. In Section \ref{sec:controlc} we \revision{divide the agents in different populations and we} introduce the control variable for the evolution of the contacts and the optimal control problem linked to it. The objective of this control variable is to minimize the relative loss of the agent's contacts. We also present an interpretation of the optimal control functional as the continuous time limit  of a discrete-time formulation. We conclude this section with the mean field description of the contacts evolution obtained through a suitable scaling of the contacts upgrade and the frequency of such evolution. Section \ref{sec:opinionmodel} is devoted to the controlled evolution of the opinions. Again, \revision{we divide the agents in distinct populations, then} we start from the microscopic multi-agent description of the binary interactions that shape each agent's opinion and present the optimal control problem linked to them. As we did for the evolution of the contacts, also for the opinion evolution we present the mean-field limit, leading to the Fokker-Planck-type equation for the time evolution of the joint density of opinions and contacts in the  regime of small updates of the opinions and the contacts happening with a very high frequency. Finally, Section \ref{sec:numerics} presents a set of numerical experiments carried on using a Monte Carlo approach, that show the different outcomes of the joint contacts and opinion evolution under various possible configurations.

\section{The evolutionary model for social network contacts}
\label{sec:contactmodel}
In this section we introduce the kinetic model governing the evolution of social media contacts. We first assume that agents’ opinions are fixed, and we specify the microscopic interaction rules governing the evolution of each agent’s number of connections. Following the approaches proposed in \cite{Albi24,Dim20}, these rules lead to a Boltzmann-type master equation describing the evolution of the corresponding distribution.

Let us consider a population of indistinguishable interacting agents, each of them characterized, at each time $t\geq 0$, by its opinion $v\in[-1,1]$ on a given topic, and by its contacts, which is the number of followers the agent has on a social media platform, that we will refer to as $c\geq 0$. 
Hence we consider the probability density of agents as $f=f(v,c,t)$ over the space of opinions and contacts at time $t>0$, such that
\[
\iint_{[-1,1]\times\mathbb{R}_+} f(v,c,t)  \mathrm{d}v \mathrm{d}c = 1, 
\]
and density marginals
\begin{equation}\label{marginals}
h(c,t):=\int_{[-1,1]} f(v,c,t) \mathrm{d}\revision{v}, \qquad g(v,t):=\int_{\mathbb{R}_+} f(v,c,t) \mathrm{d}\revision{c}.
\end{equation}
 We first analyze the evolution of agents' \revision{connections}, considering the opinions as a fixed parameter. Thus, 
for any given measurable subset $D\subseteq \mathbb{R}_+$, the quantity
\[
\int_D h(c,t)\mathrm{d}c
\]
represents the fraction of  agents whose number of contacts satisfies $c\in D$ at time $t$.

As in \cite{Albi24}, we assume that the agents are generally inclined to increase their number of social media followers once they enter the network, and that they tend to reach, at least, a certain number of followers $\overline c$, which represents a desirable popularity level. In addition, we assume that an agent’s ability to acquire or retain contacts is influenced by the popularity of the opinion they express, so that opinion–driven visibility plays a direct role in the dynamics of the contact variable.

\begin{figure}
    \centering
    \includegraphics[width=0.7\linewidth]{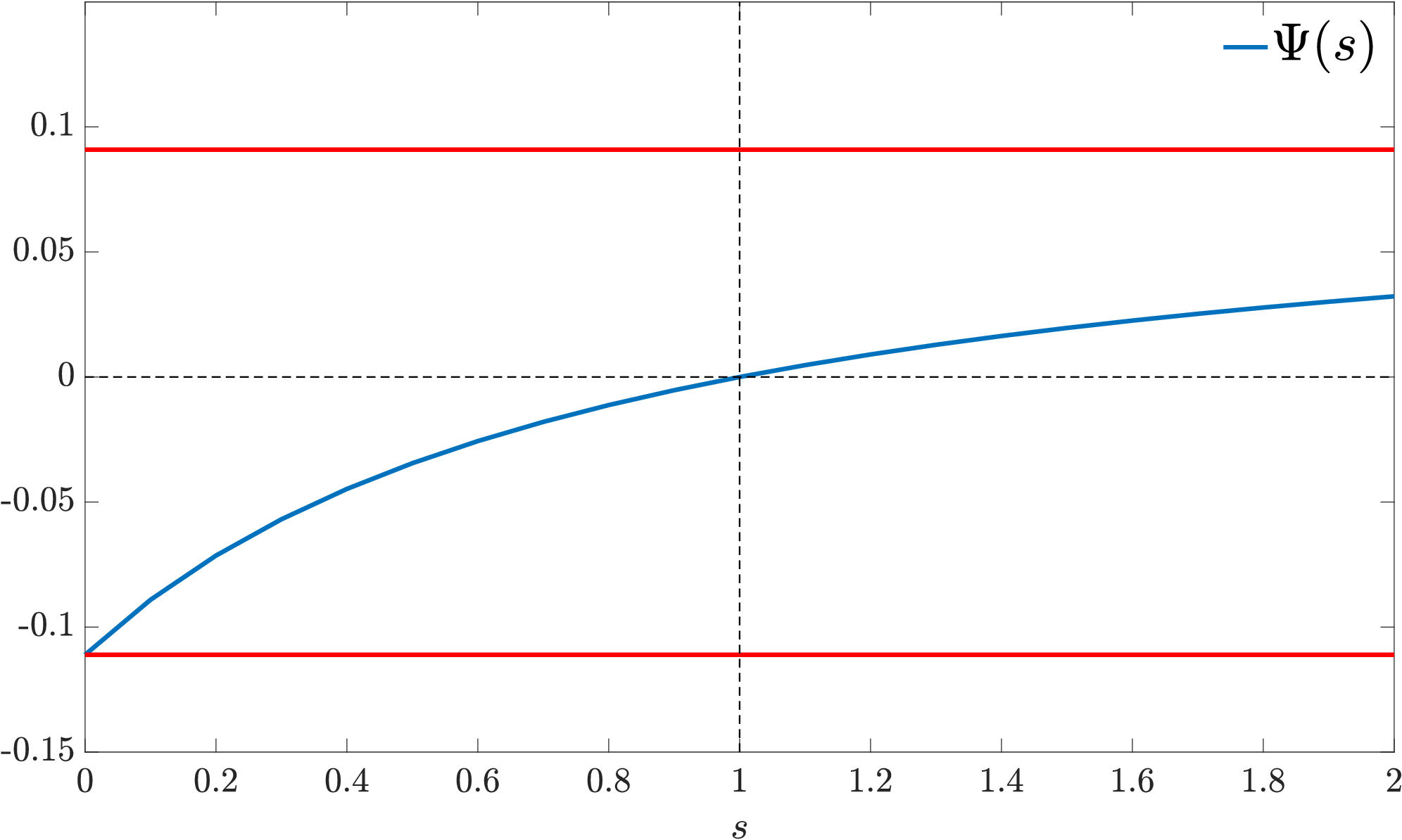}
    \caption{Plot of $\Psi(s)$. The red lines represent the bounds of the function when $\mu = 0.1$.}
    \label{fig:psi}
\end{figure}

In practice, this leads to the following microscopic dynamics for the number of connections $c$ of an agent with opinion $v$,
\begin{equation}\label{eq:contactsev}
	\revision{c^\prime = c - \left(\Psi\left(\frac{c}{\overline c}\right) + \beta\Phi(v)\right)c + \eta c},
	\end{equation}
where $\eta$ is a random variable with zero mean  and a  finite standard deviation $\nu$. Equation \eqref{eq:contactsev} models the fact that the evolution of the number of contacts is driven  by three factors. The first one is  deterministic and described by a value function $\Psi(s)$ that we suppose defined as
\begin{equation*}\label{eq:valuef}
	\revision{\Psi(s) = \left( \frac{\mu}{1-\mu}\right) \frac{s -1}{\frac{1+\mu}{1-\mu}s+1}},
\end{equation*}
with $s = c/\overline c$. This function is dimensionless, increasing and it is equal to zero at the reference point $s=1$. Moreover, it is asymmetric: moving toward the reference level $s=1$ is easier for agents starting below the threshold, while agents above the reference level face a weaker incentive to reduce their connections. This reflects the empirical principle that it is easier to gain new contacts when starting from a low number, whereas highly connected users are not motivated to decrease their popularity. This mechanism encodes the principles of the prospect theory formulated by Kahneman and Tversky \cite{KTa}. 
The function is bounded by the values
\[
\revision{-\frac{\mu}{1-\mu}\leq \Psi(s) \leq \frac{\mu}{1+\mu}},
\]
and it can be seen represented in Figure \ref{fig:psi}. This specific choice is justified by the results presented in \cite{Albi24}. 

The second contribution is also deterministic, it depends on the agent's opinion and it is represented by the function
\begin{equation*}
	\Phi (v) = \theta \left( (v-m_v)^2-\delta_\Phi^2 \right),
\end{equation*}
where $\theta>0$,  $m_v$ is the mean opinion of the entire population and $\delta_\Phi$ is a fixed parameter that depends on the tolerance of the population to unpopular opinions. This function is negative if the agent has an opinion with a distance which is less or equal than $\delta_\Phi$ from the average opinion of the population, positive \revision{otherwise}, and due to the compactness of the domain it is bounded in 
\[
-\theta \delta_{\Phi}^2 \leq \Phi(v) \leq \theta\left(4-\delta_{\Phi}^2\right).
\]

Its role in \eqref{eq:contactsev} is to decrease the popularity of the ones who do not conform with the majority of the population, while increasing the number of connections of the ones with opinion close to the mean stance of the rest of the agents.

Finally, the third contribution to the evolution of the number of connections is a stochastic one, described by the random variable $\eta$, and modeling additional sources of influence. Since the domain of the variable $c$ is $\mathbb{R}_+$, one must ensure that $\eta$ is bounded from below, in particular that
\[
\eta \geq \frac{\beta\theta\left(4-\delta_{\Phi}^2\right) \left(1+\mu\right) - 1}{1+\mu}.
\]

\subsection{The controlled evolution of social media contacts}\label{sec:controlc}

\revision{In the following, we assume that the agents are divided into a finite number of distinct populations, each of them possibly characterized by a different control strategy acting on the evolution of contacts and/or opinions.
We denote by $f_\ell(v,c,t)$ the density of agents belonging to the $\ell$-th population with opinion $v \in [-1,1]$ and number of contacts $c \in \mathbb{R}_+$ at time $t$. The total density is then given by
\[
f(v,c,t) := \sum_{\ell} f_\ell(v,c,t),
\] and the marginal densities of opinions and contacts are defined analogously to \eqref{marginals}. Accordingly, the mass of each population is defined as
\[
\rho_\ell := \iint_{[-1,1]\times \mathbb{R}_+}f_\ell(v,c,t)\, dv\, dc,\, \textrm{ with } \sum_{\ell} \rho_\ell = 1,
\]
where each population is conserved in time, since there is no exchange of agents between different populations, which can be interpreted as a compartment characterized by its own control strategy. The mean contacts and opinions write respectively as follows
\[
\begin{aligned}
\rho_\ell m_{\ell,c}(t) &:= \iint_{[-1,1]\times \mathbb{R}_+}cf_\ell(v,c,t)\, dv\, dc,\cr
\rho_\ell m_{\ell,v}(t) &:= \iint_{[-1,1]\times \mathbb{R}_+}vf_\ell(v,c,t)\, dv\, dc.
\end{aligned}
\]
}

\revision{We now describe the evolution of the agent's popularity of population $\ell$, modifying the uncontrolled dynamics  in \eqref{eq:contactsev} introducing a control input $\kappa_\ell$ in such dynamics:}
    \begin{equation}\label{eq:contactsevcontr}
		\revision{c^\prime = c -\left(\Psi\left(\frac{c}{\overline c} \right) +\beta\left(\Phi\left( v \right) - \kappa_\ell\right)\right)c + \eta c.}
	\end{equation}

\revision{In this work, we restrict ourselves to considering a common functional form of the contact control for all populations, and we distinguish the different groups only through the values of the corresponding control parameters. This choice allows us to model heterogeneous intervention capabilities within a unified control framework, while keeping the analytical structure of the model unchanged across populations.}

The optimal control $\kappa_\ell^\star$ is obtained as the minimizer of the following cost functional
\begin{equation}\label{eq:functionalc}
\kappa_\ell^\star = \arg\min_{\kappa} \, \mathcal{J}_\ell[\kappa], 
\qquad 
\mathcal{J}_\ell[\kappa] := 
\mathbb{E}\!\left[
- \lambda^\ell\,\frac{c'-c}{c}\, R^\ell_c(c)\,H^\ell_c(v)
+ \frac{\gamma_c^\ell}{2}\,\revision{\beta}\kappa^{2}
\right],
\end{equation} which penalizes the relative loss of the agent's contacts along with a quadratic regularization term. Indeed, the quantity ${(c^\prime - c)}/{c}$ is the incremental relative variation of contacts between two consecutive time steps, and overall the first term of \( \mathcal{J}_\ell \) promotes control actions capable of counterbalancing an expected loss of contacts, but only in the event that certain activation conditions are met. Specifically, such conditions are expressed by the activation functions $R^\ell_c(c)$ and $H^\ell_c(v)$. The parameter $\lambda^\ell$ sets the intrinsic scale of the control action, determining the overall intensity with which contact growth can be enhanced relative to the natural dynamics. The second term appearing in \eqref{eq:functionalc}, penalizes the control magnitude, with the positive constant $\gamma^\ell_c$ determining how costly it is to increase the number of contacts. The activation function $R^\ell_c(c)$ only depends on the current number of contacts of the agent and is modeled by the sigmoid
\begin{equation*}\label{eq:R}
	R^\ell_c(c)= \frac{1}{1+\exp(-\alpha^\ell_R(c^\ell_{\min}-c))},
\end{equation*}
which is close to one when contacts fall below the threshold \(c^\ell_{\min}\), and  decreases as \(c\) grows beyond it. The activation function $H^\ell_c(v)$ depends instead on the local opinion density $\revision{\rho}(v)$, it ensures that the control only acts on leaders positioned within a prescribed opinion region and is given by
\begin{equation*}\label{eq:H}
	H^\ell_c(v) = \frac{1}{1 + \exp(-\alpha^\ell_H(\rho(v) - \rho^\ell_\star))},
\end{equation*}
which is increasing when the agent is surrounded by many others with similar opinions and models the phenomenon of social reinforcement. With $\rho(v)$ we indicate the \emph{local opinion mass}, meaning
\[
\rho(v) =\iint_{\lvert w-v\rvert\leq r} f(w,c,t) \mathrm{d}w\mathrm{d} c ,
\]
i.e.~the fraction of agents whose opinions have a distance from $v$ which is less or equal than $r>0$, while $\rho^\ell_\star$ indicates an activation threshold, meaning the minimum mass required for a significant group effect to be witnessed. The parameters $\alpha^\ell_R,\alpha^\ell_H>0$ modify the steepness on the sigmoid functions. In this framework, $H^\ell_c(v)$ increases with the density  of similar opinions, modeling the social tendency of the agents to seek more contacts when they feel part of a group with a shared vision. 

The control model addresses the tendency to reduce the relative loss of popularity by acting when the number of agents' contacts decreases, while normalizing the action relative to the agent's current number of connections. Hence, the functional rewards any positive deviation of the contact dynamics from their uncontrolled evolution, that is, a positive relative increment of contacts.  This enhancement is activated only for agents whose popularity level  is below the desired threshold and whose opinion lies within the region where the intervention is intended to be active.

In order to obtain an explicit expression of the optimal control for a single agent, we follow a greedy approach as in \cite{Albi16}, looking for a low complexity computational representation of the optimal control problem expressed in \eqref{eq:contactsevcontr}-\eqref{eq:functionalc} using the so-called model predictive control, that computes a suboptimal solution by iterating on a sequence of finite time steps. Here, we limit ourselves to a single evolution step, from $c$ to $c^\prime$, because we are interested in an instantaneous control \revision{strategy}. For each agent, we can obtain an explicit expression of the optimal control deriving the functional $\mathcal{J}_\ell [\kappa]$ with respect to $\kappa$ and setting the derivative equal to zero. \revision{Indeed, noticing that the expected value of $\eta$ is equal to zero, and that} $\partial_\kappa c^\prime = \beta c$, we get
\[ 
\partial_\kappa\mathcal{J}_\ell[\kappa]=\mathbb{E}\!\left[
-\lambda^\ell \beta  R^\ell_c(c)H^\ell_c(v)
+ \gamma_c^\ell \revision{\beta}\kappa
\right].
\]
Imposing the first order optimality conditions $\partial_\kappa\mathcal{J}_\ell[\kappa] = 0$ and using that $\gamma^\ell_c> 0$ we get
\[
-\lambda^\ell  R^\ell_c(c)H^\ell_c(v)
+ \gamma^\ell_c \kappa =0.
\]
Finally, solving for \(\kappa\) yields the feedback control law
\begin{equation}\label{eq:controllo}
	\kappa^\star_\ell = \frac{\lambda^\ell}{\gamma_c^\ell}\,R^\ell_c(c)\,H^\ell_c(v).
\end{equation}

\begin{rmk}[Time discrete functional and continuous limit.]
We notice that the functional \eqref{eq:functionalc} can be interpreted as a one-step discrete-time functional. Indeed,
introducing a partition of the time interval \([0,T]\), given by \(t_n=n\Delta t\), \(n=0,\dots,N\), \(N\Delta t=T\), with \(\Delta t>0\), we consider analogously with \eqref{eq:functionalc} the following discrete cost functional (for the sake of notational simplicity, we omit here the population index $\ell$):
\[
\mathcal{J}_{\Delta t} \;=\; \sum_{n=0}^{N-1} \Delta t \; \mathcal{L}_n,
\qquad
\mathcal{L}_n \;=\; -\lambda\,\frac{c^{n+1}-c^n}{\Delta t\, c^n}\,R_c(c^n)H_c(v^n)
\;+\; \frac{\gamma_c \beta}{2\Delta t}\, (\kappa^n)^2 .
\]
The factor \(1/\Delta t\) inside the first term is introduced so that the first summand approximates the instantaneous relative growth rate as \(\Delta t\to 0^+\) as follows
\[
\frac{c^{n+1}-c^n}{\Delta t\, c^n} \;\approx\; \frac{\dot c(t_n)}{c(t_n)} .
\]

Hence, letting \(\Delta t\to0^+\) we obtain the continuous functional
\begin{equation}\label{eq:J_continuous}
\mathcal{J}[\kappa] \;=\; \int_0^T \Big[ -\lambda\,\frac{\dot c(t)}{c(t)}\,R_c(c(t))H_c(v(t))
\;+\; \frac{\tilde\gamma_c \beta}{2}\,\kappa(t)^2 \Big]\,\mathrm{d}t .
\end{equation} where  $\tilde\gamma_c$ is obtained by absorbing the $\Delta t$ scaling of the discrete parameter $\gamma_c$.

Observe that the integrand contains the quantity \(\dot c/c\). Formally we can write,
\[
\int_0^T \frac{\dot c(t)}{c(t)}\,\mathrm{d}t \;=\; \int_0^T \mathrm{d}\big(\ln c(t)\big)
\;=\; \ln c(T)-\ln c(0),
\]
provided \(c(t)>0\) for all \(t\) and that the functions are regular enough. Therefore, if \(R_c\) and \(H_c\) were constant in time, or independent of \(c\) and \(v\), the first integral in \eqref{eq:J_continuous} would reduce to a boundary term:
\[
- \lambda R_c H_c \int_0^T \frac{\dot c(t)}{c(t)}\,\mathrm{d}t
= -\lambda R_c H_c \big(\ln c(T)-\ln c(0)\big).
\]
This shows that maximizing the integral of \(\dot c/c\) is equivalent, up to a multiplicative constant and boundary terms, to maximizing the logarithmic growth of contacts over the time window \([0,T]\).
\end{rmk}

\subsubsection*{The Fokker-Planck type equation for the kinetic evolution of contacts}
In order to pass to a mean\revision{-}field description of contacts, we suppose that the process of network formation is the result of small increments of the number of connections of the agents. This translates in the scaling of equation \eqref{eq:contactsevcontr}\revision{.} Let $\varepsilon >0$, the scaled equation now reads
\begin{equation*}
\revision{c^\prime = c - \left(\Psi^\varepsilon \left(\frac{c}{\overline c}\right) + \varepsilon \beta\left(  \Phi(v) - \kappa_\ell\right)\right)c + \eta_\varepsilon c
}\end{equation*}
where $\eta_\varepsilon$ is a random variable with mean equal to zero,  standard deviation $\sqrt{\varepsilon}\nu$\revision{, and $\mathbb{E}\left[\left(\eta_\varepsilon\right)^3 \right] = O(\varepsilon^{3/2})$}, while the value function $\Psi^\varepsilon$ is defined as
\[
\revision{\Psi^\varepsilon(s) = \left( \frac{\mu}{1-\mu}\right) \frac{s^\varepsilon -1}{\frac{1+\mu}{1-\mu}s^\varepsilon+1}},
\]
and the scaled control is the same as in \eqref{eq:controllo}, since we both scale $\beta$ and $\gamma_c$ with the same variable $\varepsilon$.

Supposing the opinion variable $v$ fixed, we can now derive an evolutionary model for the controlled formation of the network \revision{conditioned on a given opinion.

To this end, we introduce the scaled conditional density of
contacts defined through the decomposition
$f_{\ell,\varepsilon}(v,c,t) = h_{\ell,\varepsilon}(c,t|v)g_{\ell,\varepsilon}(v,t),$
where $ g_{\ell,\varepsilon}(v,t)$ is the scaled opinion marginal and
$h_{\ell,\varepsilon}(c,t|v)$ denotes the scaled conditional distribution of
contacts among agents having opinion $v$.

In this conditional description the opinion variable acts as a ``frozen'' parameter during the contact interaction step, while the kinetic
evolution is performed only with respect to the contact variable $c$.
}
\revision{
Since the contact interaction does not modify the opinion variable, the opinion marginal $g_{\ell,\varepsilon}(v,t)$ remains unchanged during this interaction step, allowing us to write a closed kinetic equation for the conditional density $h_{\ell,\varepsilon}(c,t|v)$, which} obeys a Boltzmann-like equation of the form
\begin{equation}\label{eq:boltzc}
\frac{\mathrm{d}}{\mathrm{d}t}\int_{\mathbb{R}_+} \varphi(c)\revision{h_{\ell,\varepsilon}(c,t|v)}\mathrm{d}c = \mathbb{E}\left[\int_{\mathbb{R}_+}\chi(c)\left(\varphi(c^\prime)-\varphi(c)\right)\revision{h_{\ell,\varepsilon}(c,t|v)}\mathrm{d}c\right],
\end{equation}
where $\varphi$ is a test function and $\chi$ is the frequency of interaction, which we suppose to be very high, specifically we set it to be equal to $\chi(c) = 1/\varepsilon$. Notice that the choice $\varphi(c)=1$ leads to the conservation of the total mass of the agents, while other choices give us other information on the network evolution, for example setting $\varphi(c)=c$ leads to the evolution of the mean number of connections which we can see is not conserved in time. Assuming $\varphi$ smooth enough we can use its Taylor expansions to simplify the expression in \eqref{eq:boltzc}. We have that
\begin{equation}\label{eq:incr1}
\revision{\mathbb{E}\left[c^\prime - c\right] =  -\Psi^\varepsilon \left(\frac{c}{\overline c}\right)c - \varepsilon \beta \left(\Phi(v) - \kappa_\ell\right)c},
\end{equation}
and
\begin{equation}\label{eq:incr2}
\revision{\mathbb{E}\left[(c^\prime - c)^2\right] = \varepsilon^2\left( -\frac{1}{\varepsilon}\Psi^\varepsilon \left(\frac{c}{\overline c}\right) - \beta \left(\Phi(v) - \kappa_\ell\right)\right)^2c^2 + \varepsilon \nu^2c^2},
\end{equation}
while for the cubic increments we have
\begin{equation}\label{eq:incr3}
\begin{aligned}
\revision{\mathbb{E}\left[(c^\prime - c)^3\right] =} \,\,&\revision{\varepsilon^3\left(-\frac{1}{\varepsilon}\Psi^\varepsilon \left(\frac{c}{\overline c}\right) - \beta \left(\Phi(v) - \kappa_\ell\right)\right)^3c^3 + \mathbb{E}[\eta_\varepsilon^3]c^3 }\\
&\revision{+ 3\varepsilon^2 \nu^2 \left( -\frac{1}{\varepsilon}\Psi^\varepsilon \left(\frac{c}{\overline c}\right) - \beta \left(\Phi(v) - \kappa_\ell\right)\right)c^3}.
\end{aligned}
\end{equation}
This results in
\begin{align*}
\frac{\mathrm{d}}{\mathrm{d}t}\int_{\mathbb{R}+} &\varphi(c)\revision{h_{\ell,\varepsilon}(c,t|v)}\mathrm{d}c = \int_{\mathbb{R}_+}\left(\varphi^\prime(c) \revision{\left( -\frac{1}{\varepsilon}\Psi^\varepsilon \left(\frac{c}{\overline c}\right) - \beta\left(\Phi(v) - \kappa_\ell\right)\right)c}\right.  \left.+ \frac{\nu^2c^2}{2}\varphi^{\prime\prime}(c)\right)\revision{h_{\ell,\varepsilon}(c,t|v)}\mathrm{d}c\\
&+\revision{\frac{\varepsilon}{2}}\int_{\mathbb{R_+}}\left(\revision{\left( -\frac{1}{\varepsilon}\Psi^\varepsilon \left(\frac{c}{\overline c}\right) - \beta\left(\Phi(v) - \kappa_\ell\right)\right)^2c^2 }\varphi^{\prime\prime}(c)\right)\revision{h_{\ell,\varepsilon}(c,t|v)}\mathrm{d}c+ O(\varepsilon^{1/2}).
\end{align*}

Now letting $\varepsilon \to 0$, noticing that 
\[
\revision{\lim_{\varepsilon \to 0} \frac{1}{\varepsilon}\Psi^\varepsilon\left(\frac{c}{\overline c}\right) }= \frac{\mu}{2}\ln \left(\frac{c}{\overline c}\right)
\]
and using the notation $\revision{h_\ell(c,t|v)}$ for the limit of $\revision{h_{\ell,\varepsilon}(c,t|v)}$ as $\varepsilon \to 0$, we get the following simplified expression for the evolution of the density of the controlled connections
\begin{equation*}
\frac{\mathrm{d}}{\mathrm{d}t}\int_{\mathbb{R}+} \varphi(c)\revision{h_\ell(c,t|v)}\mathrm{d}c = \int_{\mathbb{R}_+}\left[\left(\varphi^\prime(c)\left( -\frac{\mu}{2}\ln \left(\frac{c}{\overline c}\right) - \revision{\beta\left(\Phi(v)- \kappa_\ell\right)}\right)c \right. + \frac{\nu^2c^2}{2}\varphi^{\prime\prime}(c)\right]\revision{h_\ell(c,t|v)}\mathrm{d}c,
\end{equation*}
which is the weak form of the following Fokker-Planck equation for $h_\ell(c,t|v)$
\begin{equation}\label{eq:fp}
\partial_t \revision{h_\ell(c,t|v)} = -\partial_c\left(\left( -\frac{\mu}{2}\ln \left(\frac{c}{\overline c}\right) -\revision{\beta\left(\Phi(v) - \kappa_\ell\right)}\right)c\revision{h_\ell(c,t|v)}-\frac{\nu^2}{2}\partial_{c}\left(c^2\revision{h_\ell(c,t|v)}\right)\right),
\end{equation}
\revision{coupled with zero-flux bloundary conditions in $c=0$.}

\revision{\begin{rmk}[Conditional equilibria and quasi-stationary behaviour]
In general, equation \eqref{eq:fp} does not admit an explicit global steady state due to the drift term coupling the contact dynamics with both the opinion field $\Phi(v)$ and the control $\kappa_\ell(v,c)$.

However, for any fixed $v$, \eqref{eq:fp} becomes a Fokker--Planck equation in $c$. Imposing $\partial_t h_\ell(c,t|v)=0$ and zero flux boundary conditions leads to
\[
h_{\ell,\infty}(c|v)
=
\frac{C_{\ell,\infty}(v)}{c^2}
\exp\!\left\{
-\frac{\mu}{2\nu^2}
\left(\ln\frac{c}{\overline c}\right)^2
+
\frac{2\beta}{\nu^2}
\int^c
\frac{\kappa_\ell(v,s)-\Phi(v)}{s}\,ds
\right\},
\]
with $C_{\ell,\infty}(v)$ fixed by mass conservation.

The uncontrolled mechanism produces a log-normal confinement around $\overline c$ as in \cite{Albi24}, while the control acts as an effective potential. A persistent positive drift shifts mass toward large connectivity, a negative one toward a vanishing number of contacts, and only when $\kappa_\ell(v,c)$ asymptotically balances $\Phi(v)$ then a normalizable stationary state exists.

This means that $h_{\ell,\infty}(c|v)$ represents a conditional (quasi-stationary) equilibrium: the network rapidly relaxes in $c$ given $v$, while the slow opinion dynamics continuously moves the equilibrium. A steady state arises when control and opinion feedback become independent of time and $v$.

\end{rmk}}




\section{The controlled opinion evolution}\label{sec:opinionmodel}
In this section, we focus on the evolution of opinions, represented by a continuous variable \(v \in [-1,1]\), where $-1$ indicates a strongly negative idea on a given topic, while $+1$ signifies a strongly positive one. 

\revision{Let us consider two interacting agents with pre-interaction opinions $v, v_*$ and contact levels $c, c_*$, respectively, belonging to populations $\ell$ and $k$. Notably, the two agents may belong either to the same population or to two different ones.
We denote by \(v^\prime, v_*^\prime\) the corresponding post-interaction opinions, which evolve through binary interactions depending on deterministic and stochastic factors of the form}
\begin{equation}\label{eq:binaryv} 
\begin{array}{ll}
\revision{v^\prime} &\revision{= v + \alpha \left( P(v,v_*,c,c_*)\,(v_* - v) + u_\ell \right)+ D(v)\xi,} \\
\revision{v_*^\prime} &\revision{= v_* + \alpha \left( P(v_*,v,c_*,c)\,(v- v_*) + u_{k} \right) +  D(v_*)\xi_*.}
\end{array}
\end{equation}

The random part is given by some unpredictable fluctuation in the opinion which is due, for example, to the agent's access to information through the news, and is modeled by random variables $\xi$ and $\xi_*$, with mean equal to zero, given finite standard deviation $\sigma$, \revision{and $\mathbb{E}\left[\lvert\xi\rvert^3\right]=\mathbb{E}\left[\lvert\xi_*\rvert^3\right]<+\infty$}. These variables are weighted by a self diffusion function $D(v)=1-v^2$ that is equal to zero at the boundary of the domain $[-1,1]$.
 
The deterministic part of the interaction is modeled by the compromise function \(P(v,v_*,c,c_*)\), depending on the opinions of the agents interacting and on their number of connections. The strength of this exchange is expected to be stronger when the two opinions $v$ and $v_*$ are close to each other, and when the influence associated with one individual’s connectivity exceeds that of the other. This can be mathematically modeled by expressing  the compromise function as the product of two functions
\[
P(v,v_*,c,c_*) = H(v,v_*)\,K(c,c_*),
\]
where 
\[H(v,v_*)=\mathbf{1}_{\{|v-v_*|<\delta\}}, \quad \delta >0,
\]
restricts the interactions to agents with sufficiently close opinions, and
\[\revision{
K(c,c_*) = \begin{cases}\frac{c_*^p}{c^p + c_*^p},\quad &\mbox{if } c^p+c^p_* \neq 0,\\ 0, &\mbox{otherwise},
\end{cases}}\]
weights the  relative influence of connectivity, with $p>0$. 

\revision{The other deterministic contribution is given by the control input $u_\ell$, which models an external intervention acting on the opinion of agents belonging to population $\ell$.  As for the contact control, we adopt a common functional form of the opinion control for all populations. The differences among populations are encoded through different choices of the control parameters and of the threshold functions appearing in its definition.} We denote by \revision{\(u_\ell^\star\)} the optimal control obtained as the minimizer of the quadratic cost functional defined by 
\begin{equation*}\label{eq:functionalv}
   \revision{u^\star_\ell = \arg\min_u \, \mathcal{V}_\ell[u], \qquad 
    \mathcal{V}_\ell[u] = \mathbb{E}\!\left[
        \frac{1}{2}\big(v^\prime - \tilde v_\ell \big)^2 R^\ell_v(c)H^\ell_v(v)
        + \frac{\gamma_v^\ell}{2}\alpha\, u^2
    \right]}.
\end{equation*}
The functional \revision{$\mathcal{V}_\ell$} balances two competing effects.  The first term measures the deviation of the post-interaction opinion $v^\prime$ from a prescribed target value $\tilde v_\ell$, \revision{different for each population $\ell$}, and is weighted by the activation functions \revision{$R^\ell_v(c)$ and $H^\ell_v(v)$}, which restrict the action of the control to agents that are both sufficiently active in terms of contacts and lie in a relevant opinion region.  The second term penalizes the control magnitude through the parameter $\gamma^\ell_v>0$, ensuring that excessively strong interventions are discouraged.  

Overall, the optimal feedback  control \revision{$u^\star_\ell$ aims} to steer the opinion toward the target one $\tilde v_\ell$, but only when the activation conditions encoded by \revision{$R^\ell_v(c)$ and $H^\ell_v(v)$} are satisfied. The specific choices for such  activation functions will be detailed in the numerical tests described in Section \ref{sec:numerics}. We are now interested in an explicit expression for the control $\revision{u_\ell^\star}$. First, we derive \revision{$\mathcal{V}_\ell[u]$} in $u$, noticing that \revision{$\partial_{u} v^\prime=\alpha$}, and we get
\[ 
\revision{\partial_{u}\mathcal{V}_\ell[u]=\mathbb{E}\!\left[
        \big(v^\prime - \tilde v_\ell\big) \alpha  R^\ell_v(c)H^\ell_v(v)
        + \gamma^\ell_v\alpha u
    \right]}.
\]
Then, imposing \revision{$\partial_{u}\mathcal{V}_\ell[u]=0$} and simplifying $\alpha >0$, we substitute the expression of $v^\prime$ and we pass to the expected value taking into account the fact that the random variable $\xi$ is centered in zero. In the end, we obtain
\begin{equation*}\label{eq:conto_deriv_v}
\revision{(v+\alpha(P(v,v_*,c,c_*)(v_*-v)+u) -\tilde v_\ell)R^\ell_v(c)H^\ell_v(v)+\gamma^\ell_v u = 0.}
\end{equation*}
Solving  for \(u\) yields
\begin{equation}
\label{eq:u_star}
\revision{u_\ell^\star = -\,\frac{R^\ell_v(c)H^\ell_v(v) \,\left(\,v + \alpha\left(P(v,v_*,c,c_*)(v_*-v)\right) -\tilde v_\ell\,\right)}{\gamma_v^\ell + \alpha  R^\ell_v(c)H^\ell_v(v)}.}
\end{equation}

\subsubsection*{The Fokker-Planck type equation for the kinetic evolution of opinions} 
\revision{The evolution of the density of the $\ell$-th population of agents $f_\ell(v,c,t)$ is described by the Boltzmann-like equation}
\[
\revision{\partial_t f_\ell(v,c,t) = \sum_{k} Q_{\ell}(f_\ell,f_{k})(v,c,t)},
\] 
\revision{where the Boltzmann operators $Q_{\ell}$ between populations $\ell$ and $k$, in weak form, read as}
\begin{multline}\label{eq:operQ}
\revision{\iint_{[-1,1]\times\mathbb{R}_+}\varphi(v,c)Q_{\ell}(f_\ell,f_k)\mathrm{d}v\mathrm{d}c =}\\
\revision{\iiiint_{[-1,1]^2\times \mathbb{R}^2}\mathbb{E}\left[\varphi(v^\prime,c^\prime) - \varphi(v,c)\right]f_\ell(v,c,t)f_k(v_*,c_*,t)\mathrm{d}v\mathrm{d}v_*\mathrm{d}c\mathrm{d}c_*,}
\end{multline}
\revision{with $\varphi(v,c)$ smooth test function.} Also in the opinion case, just like we did for the contacts in Section \ref{sec:controlc}, we can assume that the opinion evolution is due to small changes of the opinions happening with a high frequency. \revision{To this aim, we scale the opinion exchange intensity using the quasi-invariant scaling $\alpha\to\varepsilon\alpha$ and $\sigma^2\to\varepsilon\sigma^2$ for a small positive parameter $\varepsilon$, }
\revision{therefore the scaled version of the binary interaction \eqref{eq:binaryv} reads as}
\begin{equation*}\label{eq:binaryveps} 
\begin{cases}
v^\prime = v + \varepsilon\alpha \left( P(v,v_*,c,c_*)\,(v_* - v) + \revision{u_{\ell,\varepsilon}} \right)+ D(v)\xi^\varepsilon, \\
v_*^\prime = v_* + \varepsilon\alpha \left( P(v_*,v,c_*,c)\,(v- v_*) + \revision{u_{k,\varepsilon}} \right) +  D(v_*)\xi^\varepsilon_*,
\end{cases}
\end{equation*}
with $\xi^\varepsilon,\xi^\varepsilon_*$ random variables with mean equal to zero, variance $\mathbb{E}\left[\left(\xi^\varepsilon\right)^2\right] = \mathbb{E}\left[\left(\xi^\varepsilon_*\right)^2\right] = \varepsilon\sigma^2$\revision{, and $\mathbb{E}\left[\left(\xi^\varepsilon\right)^3\right] = \mathbb{E}\left[\left(\xi_*^\varepsilon\right)^3\right] = O(\varepsilon^{3/2})$}. 
\revision{For each population $\ell$,} the scaled density of opinion and contacts $f_{\ell,\varepsilon}(v,c,t)$ verifies the following Boltzmann-type equation
\begin{equation}\label{eq:scaledboltzv}
\revision{\partial_t f_{\ell,\varepsilon}(v,c,t) = \sum_{k} Q_{\ell}(f_{\ell,\varepsilon},f_{k,\varepsilon})(v,c,t)},
\end{equation}
\revision{where the operators $Q_{\ell}$ in \eqref{eq:scaledboltzv} are defined as in \eqref{eq:operQ}.} 

In order to write a simplified expression of \eqref{eq:scaledboltzv} we compute the Taylor expansions up to order three of both $\varphi(v^\prime,c^\prime)$ and $\varphi(v^\prime_*,c^\prime_*)$ around $v$ and $v_*$, respectively. To do so, we compute the increments
\begin{equation*}
\begin{aligned}
    \mathbb{E}\left[v^\prime - v\right] & = \varepsilon\alpha\left( P(v,v_*,c,c_*)\,(v_* - v) + \revision{u_{\ell,\varepsilon}} \right), \\
    \mathbb{E}\left[\left(v^\prime - v\right)^2\right] & = \varepsilon^2\alpha^2\left( P(v,v_*,c,c_*)\,(v_* - v) + \revision{u_{\ell,\varepsilon}} \right)^2 + \varepsilon \sigma^2 D(v)^2,
\end{aligned}
\end{equation*}
and
\[
\begin{aligned}
    \mathbb{E}\left[\left(v^\prime - v\right)^3\right] = \,\,&\varepsilon^3\alpha^3\left( P(v,v_*,c,c_*)\,(v_* - v) + \revision{u_{\ell,\varepsilon}} \right)^3 + \mathbb{E}\left[\left(D(v)\xi\revision{^\varepsilon}\right)^3\right] \\
    & + 3\varepsilon^2\alpha\sigma^2D(v)^2\left( P(v,v_*,c,c_*)\,(v_* - v) + \revision{u_{\ell,\varepsilon}} \right).
\end{aligned}
\]
Taking into consideration also \eqref{eq:incr1}, \eqref{eq:incr2}, and \eqref{eq:incr3}, and using them in \eqref{eq:scaledboltzv} \revision{, we obtain, after summing over $k$,}
\begin{equation}\label{eq:boltzvepsilon}
\begin{aligned}
    & \frac{\mathrm{d}}{\mathrm{d}t}\iint_{[-1,1]\times \mathbb{R}_+} \revision{f_{\ell,\varepsilon}}(v,c,t)\varphi(v,c)\mathrm{d}v\mathrm{d}c = \\ &\,\varepsilon\iint_{[-1,1]\times \mathbb{R}_+} \left(\alpha\partial_v\varphi(v,c)\left(\revision{\sum_{k}\mathcal{U}_{\ell,\varepsilon}[f_{k,\varepsilon}](v,c,t)}\right)\right.\\
    &\,\, + \partial_c\varphi(v,c)\left( \revision{-\frac{1}{\varepsilon}}\Psi^\varepsilon \left(\frac{c}{\overline c}\right) - \revision{\beta\left(\Phi(v) - \kappa_\ell\right)}\right)c\\
    &\,\,\,\left.+\partial_{vv}^2\varphi(v,c)\frac{\sigma^2D(v)^2}{2}+\partial_{cc}^2\varphi(v,c)\frac{\nu^2c^2}{2}\right)\revision{f_{\ell,\varepsilon}}(v,c,t)\mathrm{d}v\mathrm{d}c+ \mathcal{R}(\varepsilon),
\end{aligned}
\end{equation}
where the \revision{remainder} $\mathcal{R}(\varepsilon) = O(\varepsilon^\revision{3/2})$, and the nonlocal operator $\revision{\mathcal{U}_{\ell,\varepsilon}[f_{k,\varepsilon}]}$ is defined as
\[
\revision{\mathcal{U}_{\ell,\varepsilon}\left[f_{k,\varepsilon}\right](v,c,t) =\iint_{[-1,1]\times \mathbb{R}_+}  \left(P(v,v_*,c,c_*)(v_*-v)+u_{\ell,\varepsilon}\right)f_{k,\varepsilon} (v_*,c_*,t)\mathrm{d}v_*\mathrm{d}c_*}.
\]

\revision{We stress that the contact control $\kappa_\ell = \kappa_\ell(v,c)$ depends only on the local state variables $(v,c)$ of the tagged agent, whereas the opinion control $u_{\ell} = u_{\ell}(v,c,v_*,c_*)$ (and its rescaled version) may in general depend on the states of both interacting agents. At the kinetic level, this distinction is reflected in the limit equation: the former enters as a local drift term acting on $f_\ell$, while the latter contributes to the nonlocal drift operator $\mathcal{U}_{\ell,\varepsilon}[f_k]$ through its dependence on the interacting density.}

To model the higher frequency of interaction we can scale the time variable using the same scaling parameter $\varepsilon$, leading to the scaling $t\to\varepsilon t$. 

Now, we let $\varepsilon \to 0$. \revision{We assume that the rescaled control variable admits a finite limit denoted as $\bar u_\ell:=\lim_{\varepsilon \to 0} u_{\ell,\varepsilon}$}, and using the notation $\revision{f_\ell(v,c,t)}$ for the limit of $\revision{f_{\ell,\varepsilon}(v,c,t)}$  we have that equation \eqref{eq:boltzvepsilon} tends to
\begin{multline*}\label{eq:boltzvepsilon2}
    \frac{\mathrm{d}}{\mathrm{d}t}\iint_{[-1,1]\times \mathbb{R}_+}\revision{f_\ell(v,c,t)}\varphi(v,c)\mathrm{d}v\mathrm{d}c = \int_{\mathbb{R}_+}\int_{-1}^1\left(\alpha\partial_v\varphi(v,c)\left(\revision{\sum_{k}\mathcal{U}_\ell[f_k](v,c,t)}\right)\right.\\
     +\partial_c\varphi(v,c) \left( -\frac{\mu}{2}\ln \left(\frac{c}{\overline c}\right) -\revision{\beta\left(\Phi(v) - \kappa_\ell \right)}\right)c\\
     \left.+\partial_{vv}^2\varphi(v,c)\frac{\sigma^2D(v)^2}{2}
    +\partial_{cc}^2\varphi(v,c)\frac{\nu^2c^2}{2}\right)\revision{f_\ell(v,c,t)}\mathrm{d}v\mathrm{d}c,
\end{multline*}
with 
\begin{equation}\label{nlop}
\revision{\mathcal{U}_\ell\left[f_{k}\right](v,c,t) =\iint_{[-1,1]\times \mathbb{R}_+} \left(P(v,v_*,c,c_*)(v_*-v)+\bar u_{\ell}\right)f_{k} (v_*,c_*,t)\mathrm{d}v_*\mathrm{d}c_*},
\end{equation}
which is the weak form of the following Fokker-Planck equation
\begin{equation}\label{eq:fpv}
\begin{aligned}
     \partial_t \revision{f_\ell(v,c,t)}= &-\alpha\partial_v\left(\left(\revision{\sum_{k}\mathcal{U}_\ell[f_k](v,c,t)}\right)\revision{f_\ell(v,c,t)}\right)\\
     &-\partial_c \left( \left(-\frac{\mu}{2}\ln \left(\frac{c}{\overline c}\right) -\revision{\beta\left(\Phi(v) - \kappa_\ell\right)}\right)c\revision{f_\ell(v,c,t)}\right)\\
     &+\frac{\sigma^2}{2}\partial_{vv}^2\left(D(v)^2\revision{f_\ell(v,c,t)}\right)
    +\frac{\nu^2}{2}\partial_{cc}^2\left(c^2 \revision{f_\ell(v,c,t)}\right).
\end{aligned}
\end{equation}

\revision{\begin{rmk} [Structure of the optimal control term]In the general case, $\bar u_\ell$ may depend on both $(v,c)$ and the states of the interacting agents $(v_*,c_*)$. For the specific choice of the control \eqref{eq:u_star} we have
$$\bar u_\ell=-\,\frac{R^\ell_v(c)H^\ell_v(v) \,(\,v  -\tilde v_\ell\,)}{\gamma^\ell_v},$$ which depends only on $(v,c)$, 
and \eqref{nlop} simplifies accordingly, reducing to 
\[
\mathcal{U}_\ell[f_k](v,c,t) = \iint_{[-1,1]\times \mathbb{R}_+} P(v,v_*,c,c_*)(v_*-v) f_k(v_*,c_*,t)\, dv_* dc_* + \bar u_\ell(v,c)\rho_k.
\]
 \end{rmk}}

\revision{
\begin{rmk}[Macroscopic moment dynamics]
The evolution of macroscopic quantities is obtained from the mesoscopic model, such as \eqref{eq:fpv}, by taking velocity moments of the distribution function and integrating with respect to the microscopic variables. In this case, the presence of the nonlocal operator generally prevents the closure of the resulting moment system unless an equilibrium closure is assumed. 
Nevertheless, by imposing suitable simplifying assumptions on the nonlinear functions, one can derive a system of equations governing the first-order moments with respect to $v$. In particular, assume that $P(v,v_*,c,c_*) = 1$ and that, in the control term~\eqref{eq:u_star}, the functions $R_v(c)$ and $ H_v(v)$ are both identically equal to $1$. Then, multiplying \eqref{eq:fpv} by $v$ and integrating with respect to~$v$, we obtain the first moment equation as follows
\begin{equation}
\begin{aligned}
    \label{eq:ave_gen}
   \frac{d}{dt} m_{\ell,v}(t) &= \alpha\left(m_v(t)- m_{\ell,v}(t)\right)+\frac{\alpha}{\gamma_v^\ell}\left(\tilde v_\ell- m_{\ell,v}(t)\right)\\
    \end{aligned}
\end{equation}
resulting in a system where $m_v(t):=\sum_k\rho_km_{k,v}(t)$ denotes the average opinion of the whole population.
Furthermore, we observe that by assuming the same penalization parameter $\gamma_v^\ell \equiv \gamma_v$, multiplying by $\rho_\ell$, and summing over $\ell$ in \eqref{eq:ave_gen_all}, we can write an equation for $m_v(t)$ as follows
\begin{equation}
\begin{aligned}
    \label{eq:ave_gen_all}
   \frac{d}{dt}  m_v(t) &= \frac{\alpha}{\gamma_v}\left(\sum_\ell\rho_\ell\tilde v_\ell - m_v(t)\right),
    \end{aligned}
\end{equation}
which corresponds to a relaxation toward the weighted average of the desired target opinions $\tilde v_\ell$ of the different populations.
\end{rmk}}
\section{Numerical experiments}\label{sec:numerics}

In this Section we present a set of numerical experiments aimed at illustrating the behavior of the proposed model and analyzing the joint evolution of contacts and opinions in the presence of the control mechanisms introduced in the previous sections.

\revision{In the multi-population setting described above, the state of the system is given by a family of distribution functions $\{f_\ell(v,c,t)\}_{\ell}$, one for each population, whose evolution is governed by a system of Boltzmann-type equations~\eqref{eq:boltzc} coupled through cross-interaction operators. In the Fokker--Planck regime~\eqref{eq:fpv}, we approximate this system by means of a Monte Carlo–type strategy based on an asymptotic stochastic particle method, suitable for kinetic equations in the quasi-invariant limit.

More precisely, denoting by $f_\ell^n$ the numerical approximation of $f_\ell$ at time $t^n$, we consider the discrete update
\begin{equation}\label{eq:dsmc}
    f_{\ell}^{n+1} 
    = \left(1 - \frac{\Delta t}{\varepsilon}\right) f_{\ell}^n
      + \frac{\Delta t}{\varepsilon}\, \sum_{k} \mathcal{G}_{\ell,\varepsilon}\!\left(f_\ell^n, f_k^n\right),
\end{equation}
where the gain operators $\mathcal{G}_{\ell,\varepsilon}$ account for the inflow of particles in population $\ell$ at  $(v,c)$ due to interactions with agents of population $k$, according to the microscopic rules~\eqref{eq:contactsevcontr} and~\eqref{eq:binaryv}.
The corresponding particle scheme is detailed in Algorithm~\ref{nanbu}, where, for simplicity, we take $\varepsilon = \Delta t$. Further details on this class of asymptotic particle methods for kinetic equations in the quasi-invariant limit can be found, e.g., in~\cite{PT13}.}

\begin{algorithm}
\caption{Asymptotic particle-based scheme (Nanbu-like algorithm) \revision{with multiple populations}}\label{nanbu}
\begin{algorithmic}[h!]

\State Fix $0<\varepsilon=\Delta t<1$ and $N_s$.
\State \revision{Sample $\{v_i^0,c_i^0,\ell_i\}_{i=1}^{N_s}$ from the initial multi-population distribution $\{f^0_\ell(v,c)\}$.}

\For{$n = 0:N_t-1$}

    \State Set $N_c = \mathrm{round}(N_s/2)$.
    \State Select $N_c$ unordered random pairs $(i,i_*)$ uniformly without repetition.
    
    \For{$k = 1:N_c$}
        \State Let $(i,i_*)$ be the $k$-th selected pair.
        \State \revision{Let $\ell = \ell_i$ and $\ell_* = \ell_{i_*}$.}

        \State Sample noise terms $\eta_{\varepsilon,i}^n,\eta_{\varepsilon,i_*}^n$ with 
        $\langle\eta\rangle=0$, \ $\langle \eta^2\rangle=\varepsilon\nu^2$, and $\xi_i^n, \xi_{i_*}^n \sim \mathcal{N}(0,1)$.

        \State Compute the contact controls \revision{$\kappa_i^{n}=(\kappa_\ell)_i^{n}$ and $\kappa_{i_*}^{n}=(\kappa_{\ell_*})_{i_*}^{n}$. 
        \[
        \kappa_i^{n}
        = \frac{ \lambda^{\ell}}{\gamma^{\ell}_c} R^{\ell}_c(c^n_i)H^{\ell}_c(v^n_i),
        \qquad
        \kappa_{i_*}^{n}
        =  \frac{\lambda^{\ell_*}}{\gamma^{\ell_*}_c} R^{\ell_*}_c(c^n_{i_*})H^{\ell_*}_c(v^n_{i_*}).
        \]}

        \State Compute the opinion controls \revision{$u_i^{n}=(u_\ell)_i^{n}$ and $u_{i_*}^{n}=(u_{\ell_*})_{i_*}^{n}$.
        \begin{equation*}
        \begin{aligned}
        u_i^{n}
        &= - \frac{R_v^{\ell}(c_i^n)H_v^{\ell}(v_i^n) 
        \left(v_i^n+\varepsilon\alpha P(v_i^n,v_{i_*}^n,c_i^n,c_{i_*}^n)(v_{i_*}^n - v_i^n) - \tilde v^{\ell}\right)}
        {\gamma_v^{\ell} + \varepsilon \alpha R_v^{\ell}(c_i^n)H_v^{\ell}(v_i^n)}, \\
        u_{i_*}^{n}
        &= - \frac{R_v^{\ell_*}(c_{i_*}^n)H_v^{\ell_*}(v_{i_*}^n) 
        \left(v_{i_*}^n+\varepsilon\alpha P(v_{i_*}^n,v_i^n,c_{i_*}^n,c_i^n)(v_i^n - v_{i_*}^n) - \tilde v^{\ell_*}\right)}
        {\gamma_v^{\ell_*} + \varepsilon \alpha R_v^{\ell_*}(c_{i_*}^n)H_v^{\ell_*}(v_{i_*}^n)}.
        \end{aligned}
        \end{equation*}}
        
        \State Update contacts.
        \[
        c_j^{n+1}
        = c_j^n -\Psi^\varepsilon\!\left(\frac{c_j^n}{\overline c}\right)
        - \varepsilon\beta\!\left(\Phi(v_j^n) - \kappa_j^n\right)c_j^n
        + \eta_{\varepsilon,j}^n\, c_j^n ,
        \qquad  j\in\{i,i_*\}.
        \]

        \State Update opinions.
        \begin{equation*}
        \begin{aligned}
        v_i^{n+1} 
            &= v_i^{n} 
            + \varepsilon \alpha \!\left( P(v_i^n,v_{i_*}^n,c_i^n,c_{i_*}^n)(v_{i_*}^n - v_i^n) + u_i^n \right)
            + \sqrt{\varepsilon}\,\sigma\, D(v_i^n)\, \xi_i^n , \\[2mm]
        v_{i_*}^{n+1} 
            &= v_{i_*}^{n} 
            + \varepsilon \alpha \!\left( P(v_{i_*}^n,v_i^n,c_{i_*}^n,c_i^n)(v_{i}^n - v_{i_*}^n) + u_{i_*}^n \right)
            + \sqrt{\varepsilon}\,\sigma\, D(v_{i_*}^n)\, \xi_{i_*}^n .
        \end{aligned}
        \end{equation*}

    \EndFor

\EndFor

\end{algorithmic}
\end{algorithm}


\subsection{Test 1: comparative analysis of control strategies in coupled contact--opinion dynamics}

In this test we investigate the role of contact and opinion control mechanisms for two \revision{populations} of agents, structured as a ``leader--follower'' system, where leaders denote agents endowed with control authority (or externally driven dynamics), while followers are agents that evolve according to interaction rules and respond to the influence of the leaders. \revision{In particular, $25\%$ of the agents are designated as ``leaders'' (population $\ell=L$), whereas the others represent the ``followers'' (population  $\ell=F$).}

\begin{table}[h!]
\centering
\caption{Summary of model parameters, \revision{the superscripts in the control terms correspond to the population $\ell=L$ for leaders and $\ell=F$ for followers.}}\label{tab1}
\begin{tabular}{ll}
\toprule
\textbf{Parameter} & \textbf{Value} \\
\midrule

\multicolumn{2}{l}{\textbf{Contacts dynamics}} \\
Interaction strength & $\beta = 1$ \\
Endogenous modulation $\Psi(\cdot)$ & \revision{$\mu = 0.25$, $\bar c = 100$} \\
Endogenous modulation $\Phi(\cdot)$ & $\theta = 2$, $\delta_\Phi = 0.1$ \\
Contacts noise & $\eta$, $\nu = 0.1$ \\
\midrule

\multicolumn{2}{l}{\textbf{Opinion dynamics}} \\
Binary interaction & $\alpha = 1$, $\delta = 0.8$, $p = 3$ \\
Opinions noise & $\xi$, $\sigma = 0.1$ \\
\midrule

\multicolumn{2}{l}{\textbf{Contact control $\kappa^\star$}} \\
Cost penalization & $\gamma^L_c = 1$ \\
Maximal control intensity & $\lambda^F = 0$, $\lambda^L = 1$ (when activated) \\
Activation thresholds & $\alpha^L_R = \alpha^L_H = 0.1$, $r^L = 0.7$, $\rho_\star^L = 0.5$ \\
Contact activation level & \revision{$c^L_{\min} = 150$} \\
\midrule

\multicolumn{2}{l}{\textbf{Opinion control $u^\star$}} \\
Thresholds  & $R^F_v(c) = H^F_v(v) = 0$, $R^L_v(c) = H^L_v(v) = 1$ \\
Target opinion & $\tilde v^L = 0.5$ \\
Cost penalization & $\gamma_v^L = 10$ \\
\bottomrule
\end{tabular}
\end{table}

\begin{figure}[htbp]
\centering
\includegraphics[width=0.4\textwidth]{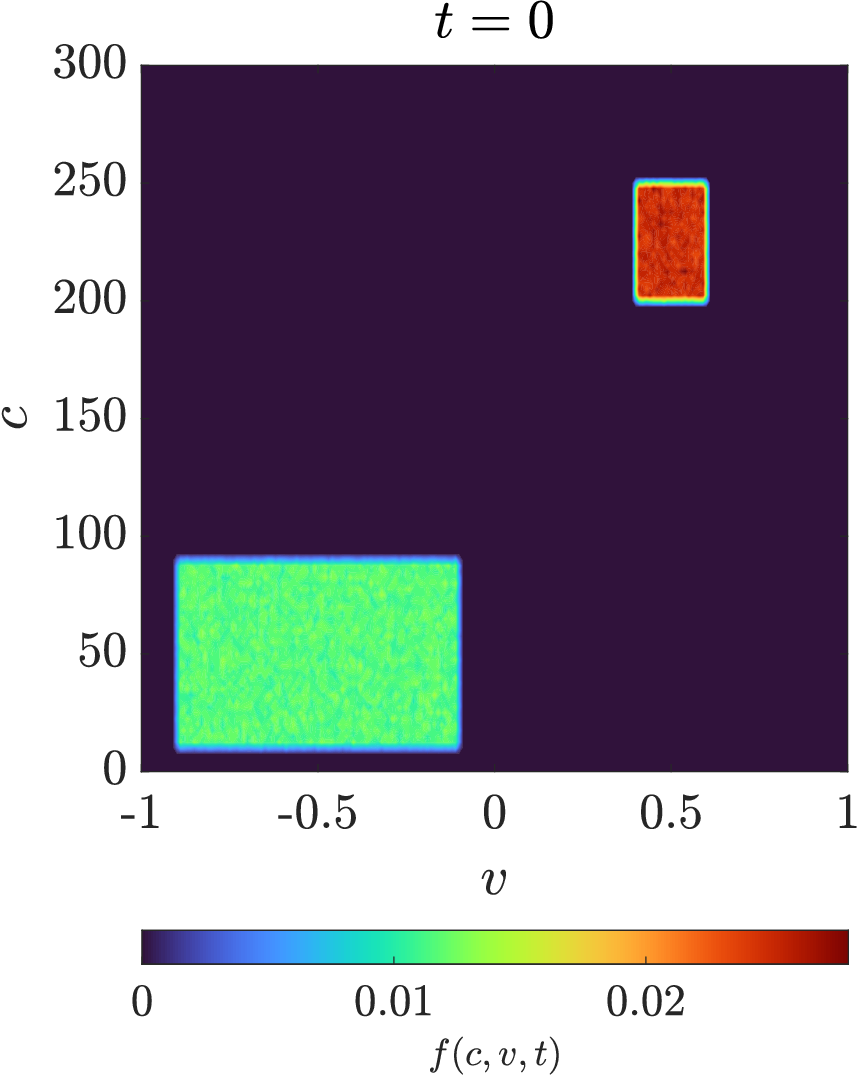}
\caption{ Test 1. Joint density $f(c,v,0)$ at the initial time.}
\label{fig:test1_t0}
\end{figure}

We sample $N_s = 10^6$ agents to simulate the model. At the initial time $t=0$, the total joint density $f_0(c,v)=f(c,v,0)$ is represented in Fig.~\ref{fig:test1_t0},  and is provided as follows: the leaders' states $(c,v)$ are sampled from a continuous uniform distribution
\(
(c,v) \sim \mathrm{Unif}\big([200,250]\times [0.4,0.6]\big),
\)
representing a cohesive and highly connected group with an opinion significantly different from that of the mass, which instead is initialized according to
\(
(c,v) \sim \mathrm{Unif}\big([10,90]\times [-0.9,-0.1]\big),
\)
corresponding to a low-connectivity population spread over a broad range of opinions. 

The simulation runs up to a final time $T = 50$ with uniform time step $\Delta t = 10^{-3}$.

In this test the followers are always uncontrolled, while for the leaders we switch on/off the contact control $\kappa^\star$ and the opinion control $u^\star$, obtaining the four scenarios (a)--(d) described below:
(a) fully uncontrolled dynamics for both populations;
(b) contact control $\kappa^\star$ active only on the leader population;
(c) opinion control $u^\star$ active only on the leader population;
(d) both controls simultaneously active on the leader population.
In all cases, the follower population evolves without control.

All numerical parameters are listed in Table~\ref{tab1}; only the activation of the control terms changes across scenarios.

\begin{figure}[htbp]
\centering

\begin{subfigure}[b]{\textwidth}
    \centering
    \includegraphics[width=\textwidth]{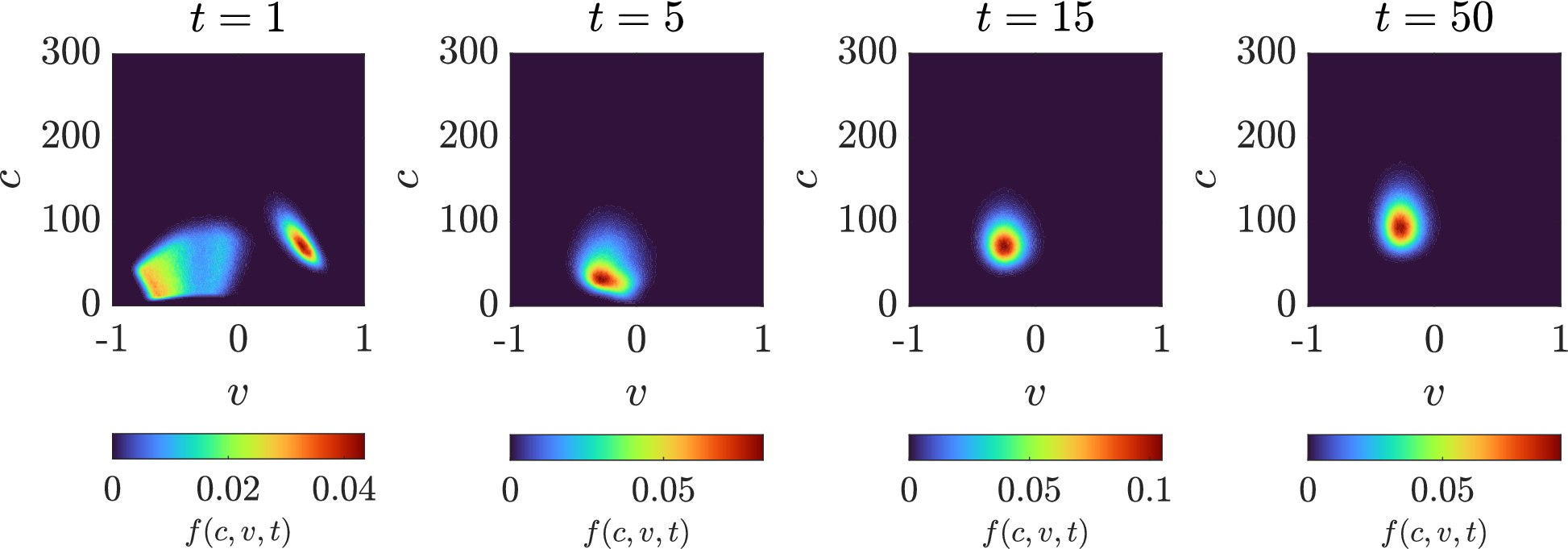}
    \caption{No control: free evolution of $f(c,v,t)$.}
    \label{fig:sim0}
\end{subfigure}

\par\vspace{0.5em}

\begin{subfigure}[b]{\textwidth}
    \centering
    \includegraphics[width=\textwidth]{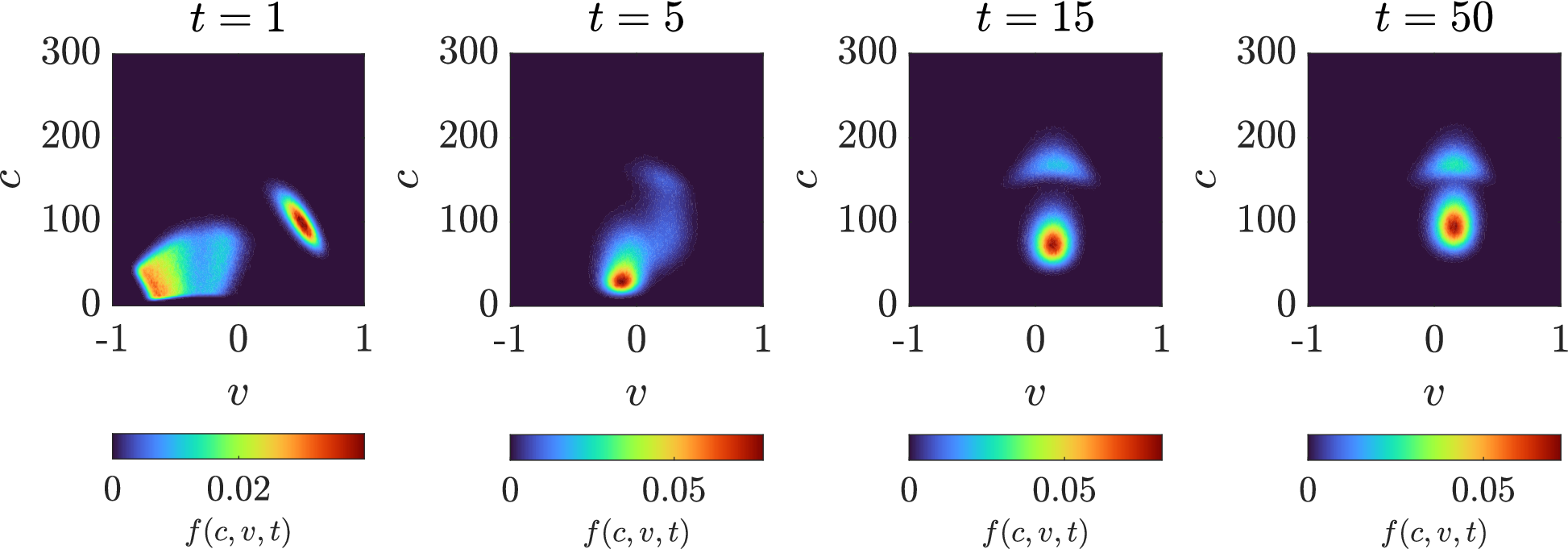}
    \caption{Control on contacts only ($\kappa^\star$ active).}
    \label{fig:sim1}
\end{subfigure}

\par\vspace{0.5em}

\begin{subfigure}[b]{\textwidth}
    \centering
    \includegraphics[width=\textwidth]{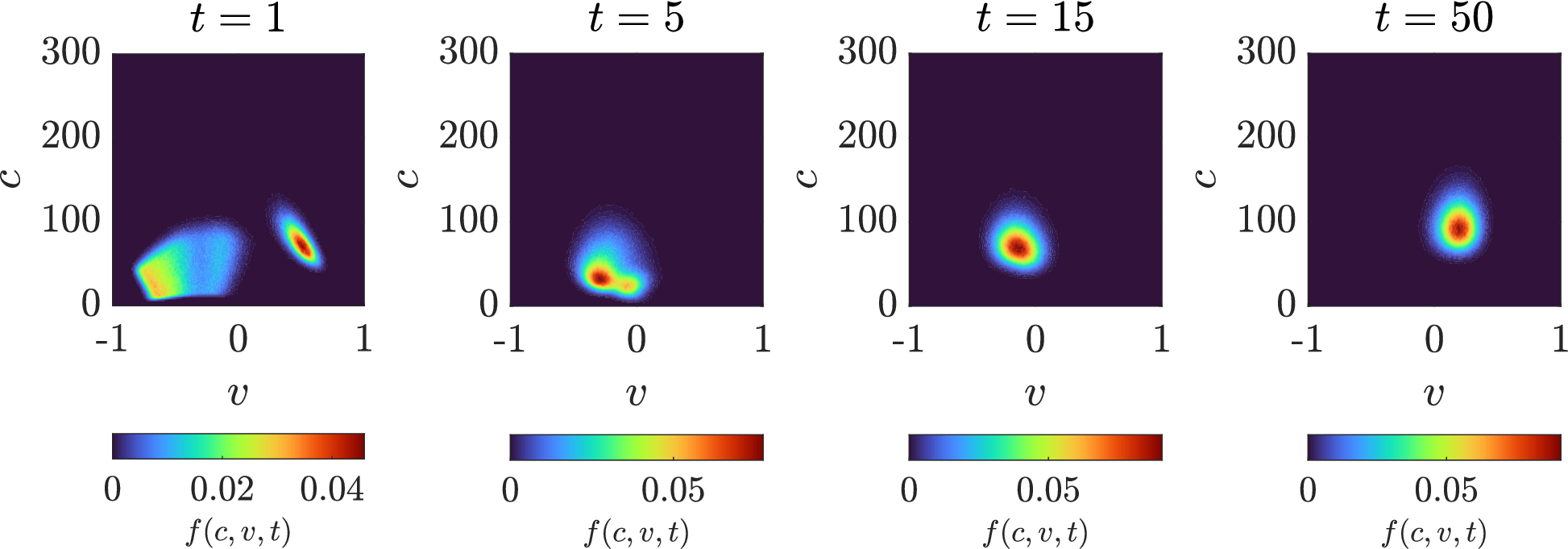}
    \caption{Control on opinions only ($u^\star$ active).}
    \label{fig:sim2}
\end{subfigure}

\par\vspace{0.5em}

\begin{subfigure}[b]{\textwidth}
    \centering
    \includegraphics[width=\textwidth]{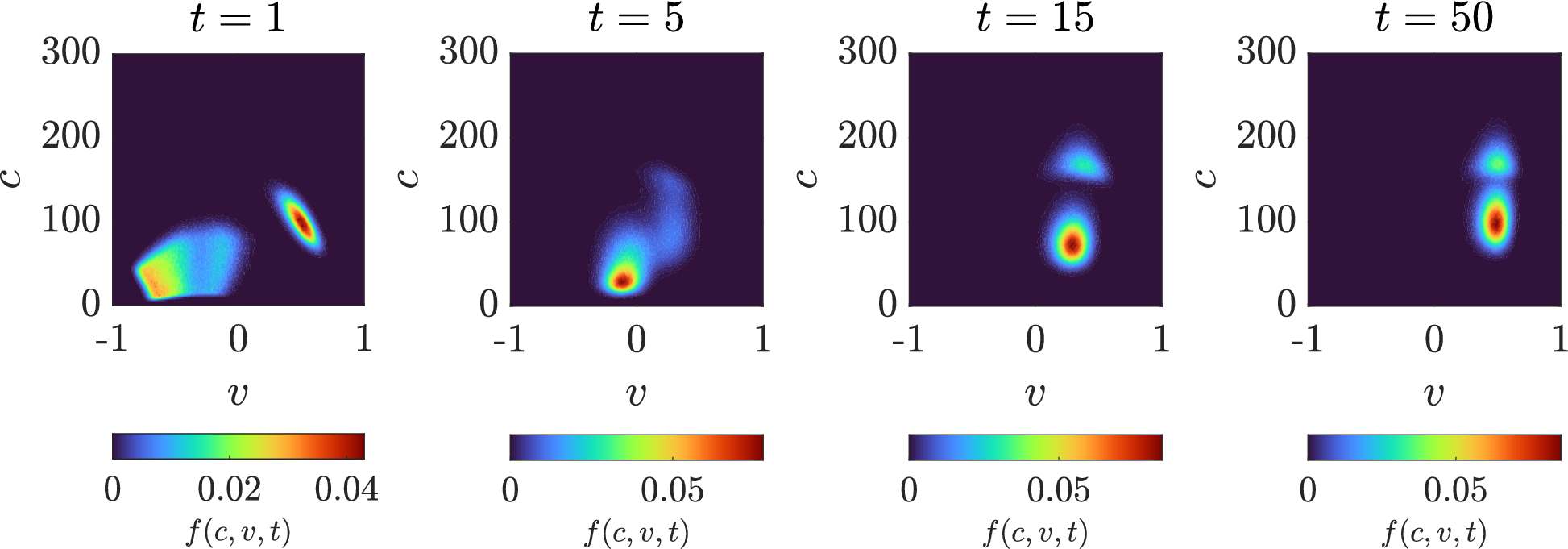}
    \caption{Joint control on contacts and opinions.}
    \label{fig:sim3}
\end{subfigure}

\caption{ Test 1. 
Time evolution of the joint density $f(c,v,t)$ in the $(v,c)$ plane for different control strategies.
Each panel shows snapshots at times $t=1, 5, \revision{15}, 50$. 
(a) uncontrolled dynamics; 
(b) control on contacts; 
(c) control on opinions; 
(d) both controls active.
}
\label{fig:simulations}
\end{figure}

\begin{figure}[htbp]
\centering

\begin{subfigure}[b]{0.48\textwidth}
    \centering
    \includegraphics[width=\textwidth]{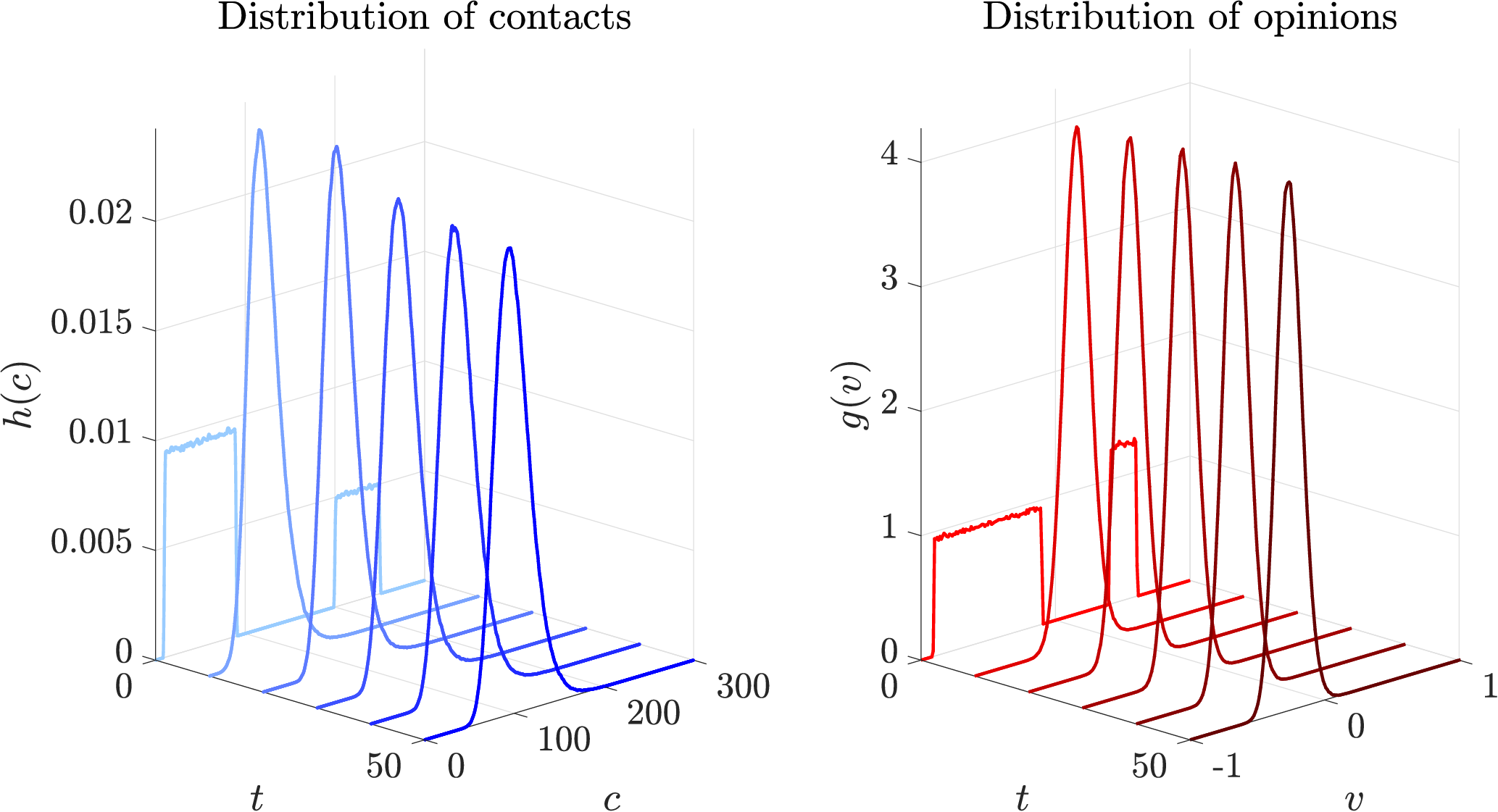}
    \caption{No control: free evolution.}
    \label{fig:dist0}
\end{subfigure}
\hfill
\begin{subfigure}[b]{0.48\textwidth}
    \centering
    \includegraphics[width=\textwidth]{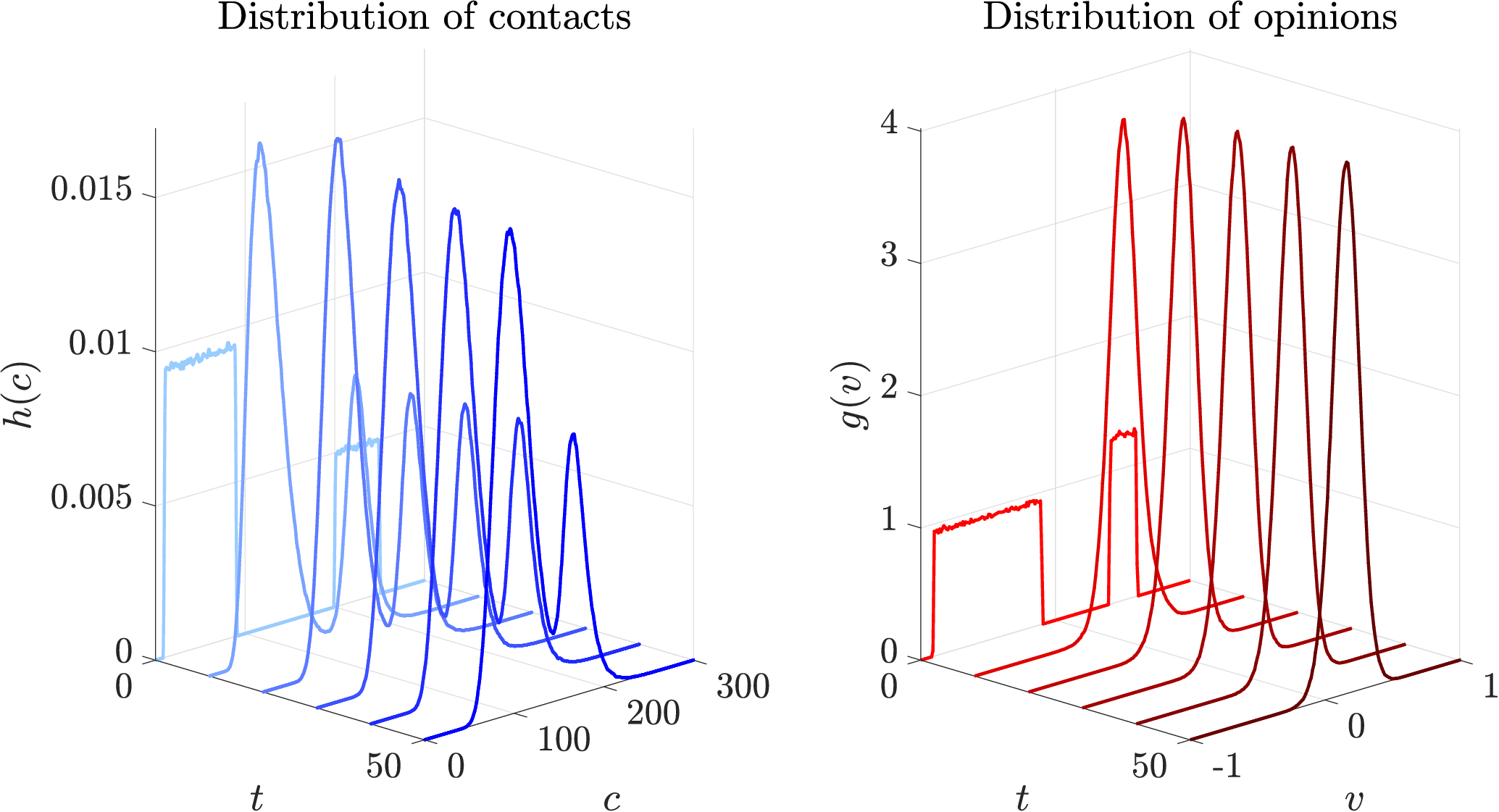}
    \caption{Control on contacts only ($\kappa^\star$ active).}
    \label{fig:dist1}
\end{subfigure}

\par\vspace{0.5em}

\begin{subfigure}[b]{0.48\textwidth}
    \centering
    \includegraphics[width=\textwidth]{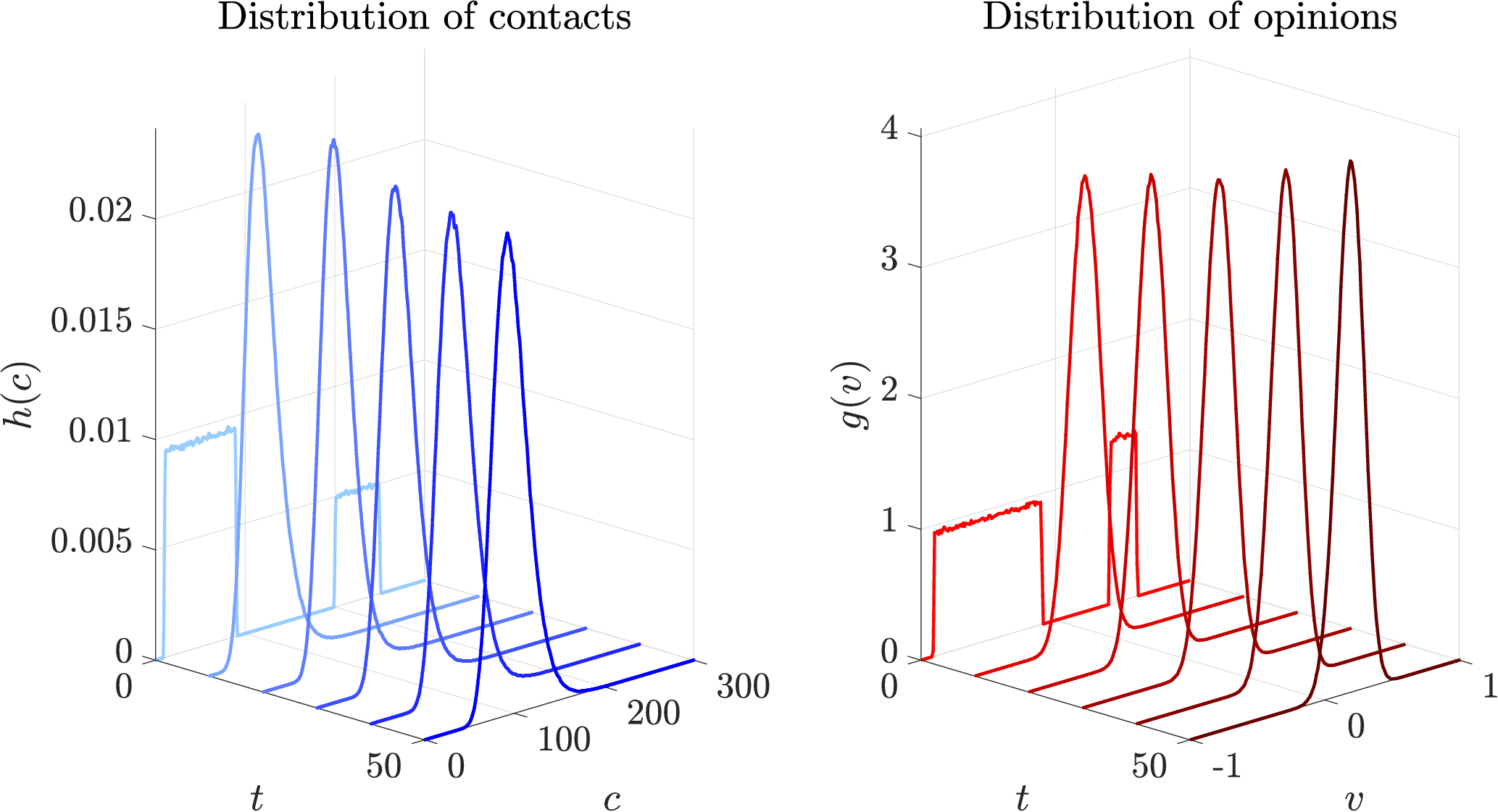}
    \caption{Control on opinions only ($u^\star$ active).}
    \label{fig:dist2}
\end{subfigure}
\hfill
\begin{subfigure}[b]{0.48\textwidth}
    \centering
    \includegraphics[width=\textwidth]{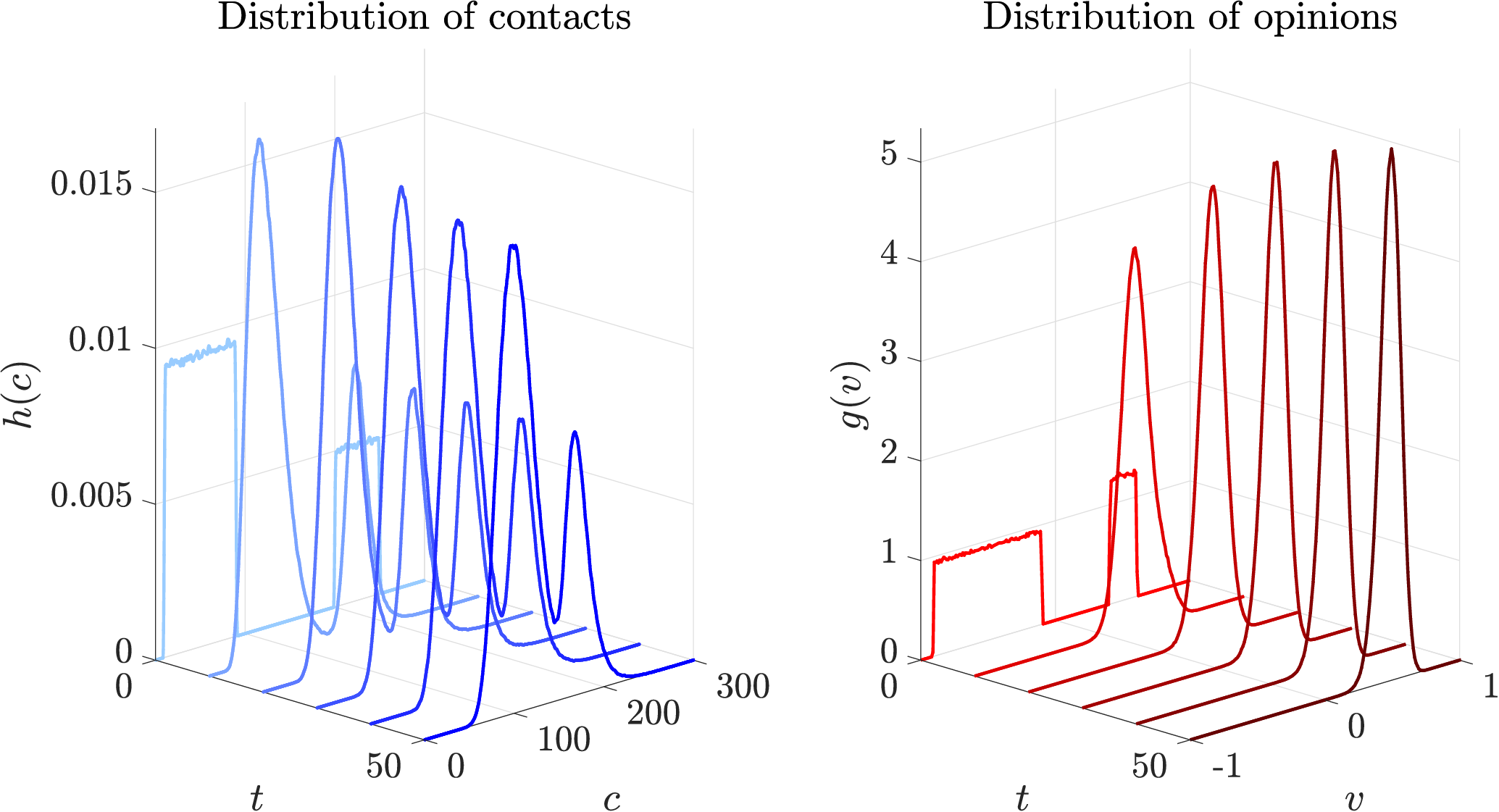}
    \caption{Joint control on contacts and opinions.}
    \label{fig:dist3}
\end{subfigure}

\caption{Test 1.
Time evolution of the marginal distributions of contacts and opinions at different times for different control strategies. 
Each panel shows both distributions together: contacts (left) and opinions (right). 
(a) uncontrolled dynamics; 
(b) control on contacts; 
(c) control on opinions; 
(d) both controls active.
}
\label{fig:distributions}
\end{figure}

\begin{figure}[htbp]
\centering
\includegraphics[width=\textwidth]{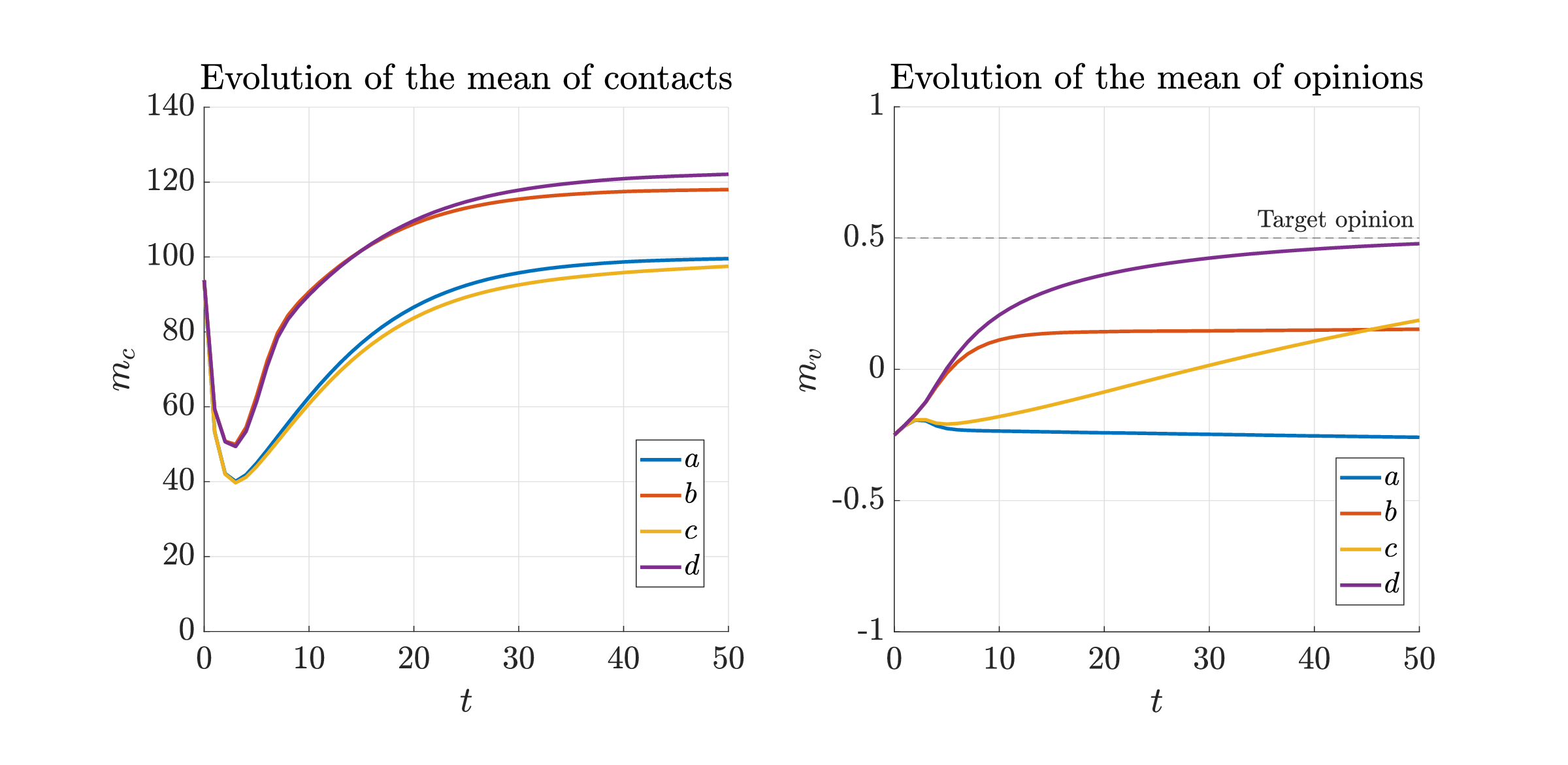}
\caption{Test 1. Time evolution of the mean number of contacts $m_c(t)$ (left panel) and the mean opinion $m_v(t)$ (right panel) over the interval $[0,T]$. In both panels we report the temporal averages corresponding to the four scenarios (a)--(d), each represented with a different colour.}

\label{fig:medie1}
\end{figure}

Simulation results are reported in Fig.~\ref{fig:simulations}, which displays snapshots of the total joint density $f(c,v,t)$ at selected times for each scenario, and in Fig.~\ref{fig:distributions}, illustrating the temporal evolution of the marginal distributions of contacts $h(c\revision{,t})$ and opinions $g(v\revision{,t})$.

In the uncontrolled case (a), the leader group, initially located around an extreme opinion, progressively loses contacts, thereby reducing its influence, and is eventually absorbed into the followers' cluster, converging to a compromise opinion.  
In scenario (b), leaders initially lose popularity, but the contact control activates near the threshold $c^L_{\min}$, preventing collapse and enabling them to regain influence; as a result, the overall opinion shifts toward a compromise closer to the leaders' initial stance.  
Scenario (c), with only opinion control active, shows that the loss of contacts dominates the dynamics: leaders rapidly lose influence and merge into the followers’ distribution, again converging to a compromise.  
Finally, in scenario (d), the combined action of both controls preserves leaders’ popularity and enables them to steer the population toward their target opinion.

In Figure~\ref{fig:medie1} we report the time evolution of the mean number of contacts $m_c(t)$ and of the mean opinion $m_v(t)$ for the four scenarios (a)--(d).  
The plots show that when the contact control is active, namely in scenarios (b) and (d), a larger and more persistent increase in the global mean of contacts is displayed, compared to the remaining cases.  
Moreover, in scenario (d), where both opinion and contact controls are simultaneously activated, the combined action of the two control layers drives the population's mean opinion towards the target value promoted by the leader group, highlighting the synergistic effect of the coupled control mechanisms.

\begin{figure}[htbp]
\centering
\includegraphics[width=\textwidth]{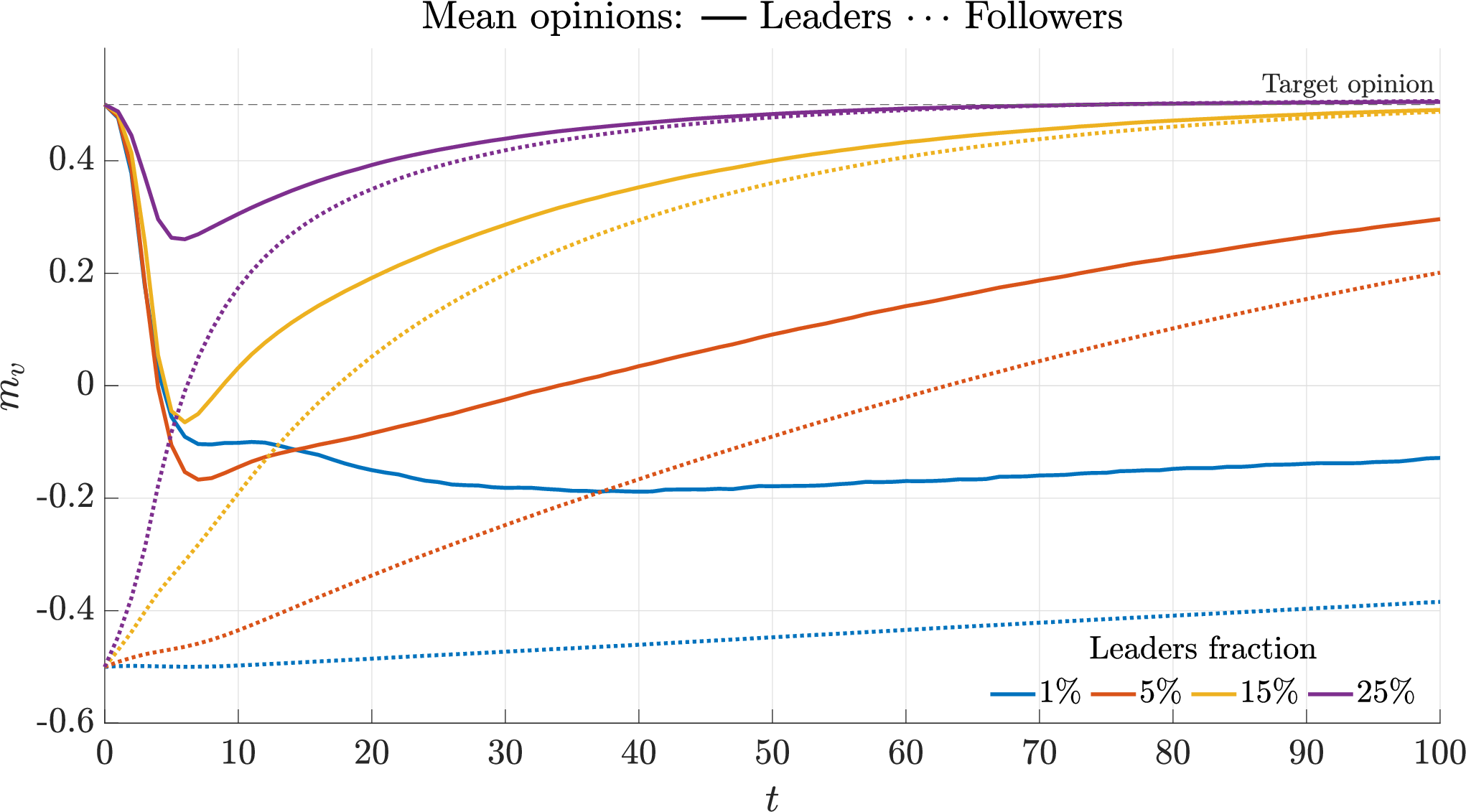}
\caption{\revision{Test 1. Time evolution of the mean opinion $m_v(t)$ over the interval $[0,100]$ for Leaders and Followers in scenario (d). Each color corresponds to a different fraction of leaders (1\%, 5\%, 15\%, 25\%). Line styles distinguish the two populations: continuous lines represent Leaders, while dotted lines represent Followers. The horizontal dashed line indicates the target opinion at $0.5$.}}
\label{fig:perc1}
\end{figure}

\revision{To address the effect of the fraction of leaders, we repeated scenario (d) varying the initial proportion of leaders from 1\% up to 25\%, in a larger temporal interval with \(T=100\).  We observe that as the fraction of leaders increases, the attraction of the population towards the target opinion promoted by the leaders becomes stronger. Interestingly, at the initial stages, the opinions of the leaders tend to slightly diverge from their own mean target opinion,  while the followers' opinions progressively align with those of the leaders, the rate of alignment depending on the strength of the attraction mechanism.  These dynamics are illustrated in Fig.~\ref{fig:perc1}.}

This test highlights how influence depends on visibility. When leaders let their popularity decay, their capacity to shape collective opinion quickly disappears, and they are reabsorbed into the majority view. Maintaining contacts prevents this loss of visibility, allowing leaders to maintain a role in directing the opinion landscape. The combination of preserving visibility and exerting persuasive effort produces lasting influence, in this case persuasion alone is ineffective when the audience becomes too small.

\subsection{Test 2: leader-driven opinion steering: symmetric vs. asymmetric competition}

In this test we study the effect of opinion control in a social system consisting of two competing leader groups, each promoting an opposite opinion, and a third uncontrolled group representing the mass of followers.  In particular, we analyze how symmetric or asymmetric penalization of the opinion-control effort influences the outcome of the competition for consensus.

\revision{We consider three populations of  agents divided as follows: $25\%$ constitute Leader Group~A ($\ell=A$), another $25\%$ constitute Leader Group~B ($\ell=B$), and the remaining $50\%$ form the ``mass'' ($\ell=M$).}  
The simulation horizon is $T=50$ with time step $\Delta t = 10^{-3}$.

The evolution of contacts and opinions follows the same parameters used in Test~1 as in Tab~\ref{tab1}.  The only difference concerns the penalization parameters for the opinion control, indicated here by $\gamma_v^A$ and $\gamma_v^B$ for the two leader populations and for the target opinions $\tilde v_A, \tilde v_B$.

\begin{figure}[htbp]
\centering
\includegraphics[width=0.4\textwidth]{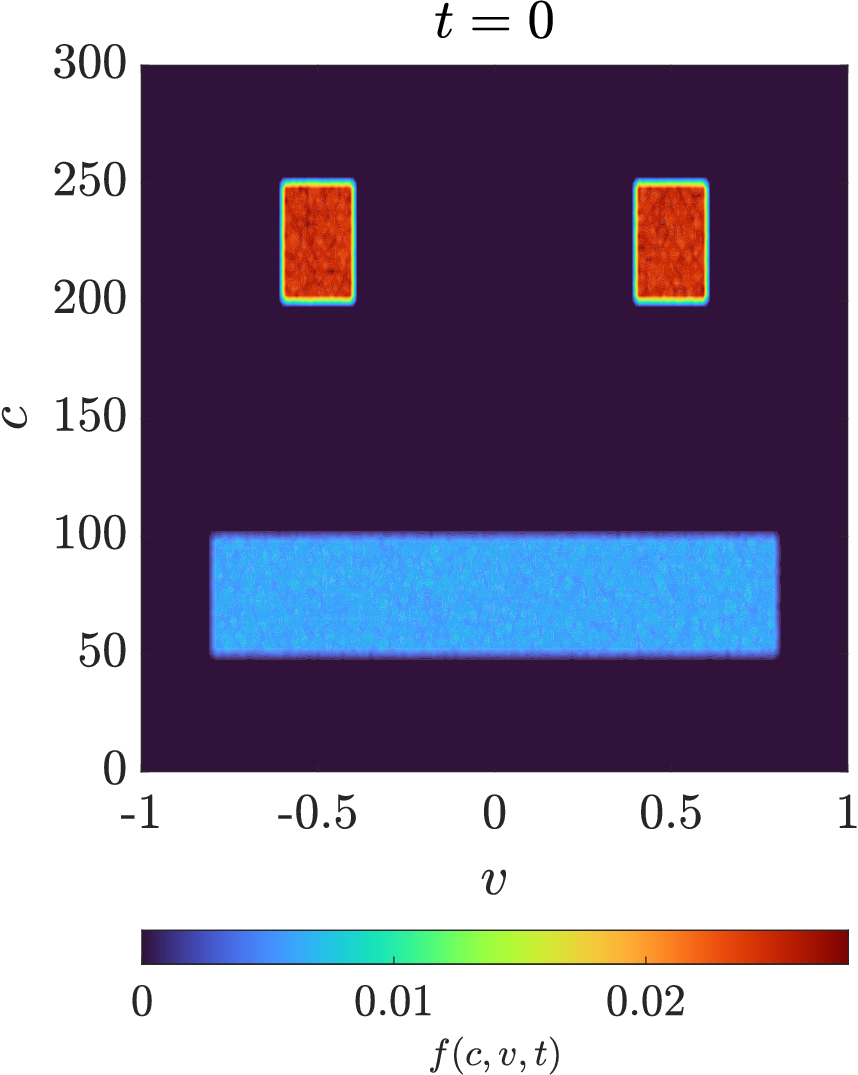}
\caption{ Test 2. Joint density $f(c,v,0)$ at the initial time.}
\label{fig:test2_t0}
\end{figure}

In Fig.~\ref{fig:test2_t0} we report the initial datum $f_0(c,v)$, sampled from $N_s = 10^6$ agents. At the initial time $t=0$, leaders of Group~A are sampled from
\(
(c,v) \sim \mathrm{Unif}\big([200,250]\times [-0.6,-0.4]\big), \)
with a target opinion  \(\tilde v_A = -0.5,
\)
representing a highly connected group promoting a moderately negative opinion.  
Leaders of Group~B are sampled from
\(
(c,v) \sim \mathrm{Unif}\big([200,250]\times [0.4,0.6]\big),\)
and similarly promote a positive target opinion \(\tilde v_B = 0.5
\).  
Followers are initialized according to
\(
(c,v) \sim \mathrm{Unif}\big([50,100]\times[-0.8,0.8]\big),
\)
thus representing a poorly connected population spread across a wide range of opinions.

We analyse three scenarios: (a) no control: both leader groups evolve freely under the uncontrolled opinion dynamics; (b) symmetric control: both leader groups employ the optimal opinion-control strategy with identical penalization, $\gamma_v^A=\gamma_v^B=1$; (c) asymmetric control: group~A faces a much higher penalization, $\gamma_v^A=100$, while group~B keeps $\gamma_v^B=1$.

\begin{figure}[htbp]
\centering

\begin{subfigure}[b]{\textwidth}
    \centering
    \includegraphics[width=\textwidth]{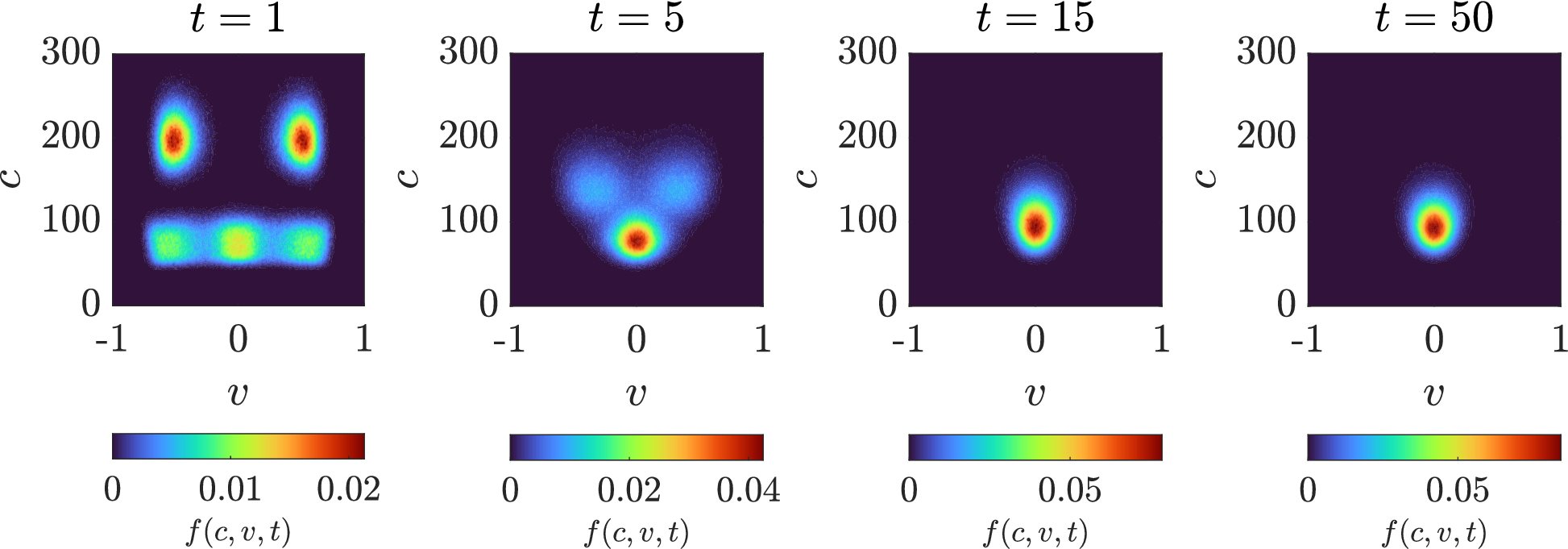}
    \caption{No control: free evolution of $f(c,v,t)$.}
    \label{fig:test2_sim0}
\end{subfigure}

\par\vspace{0.5em}

\begin{subfigure}[b]{\textwidth}
    \centering
    \includegraphics[width=\textwidth]{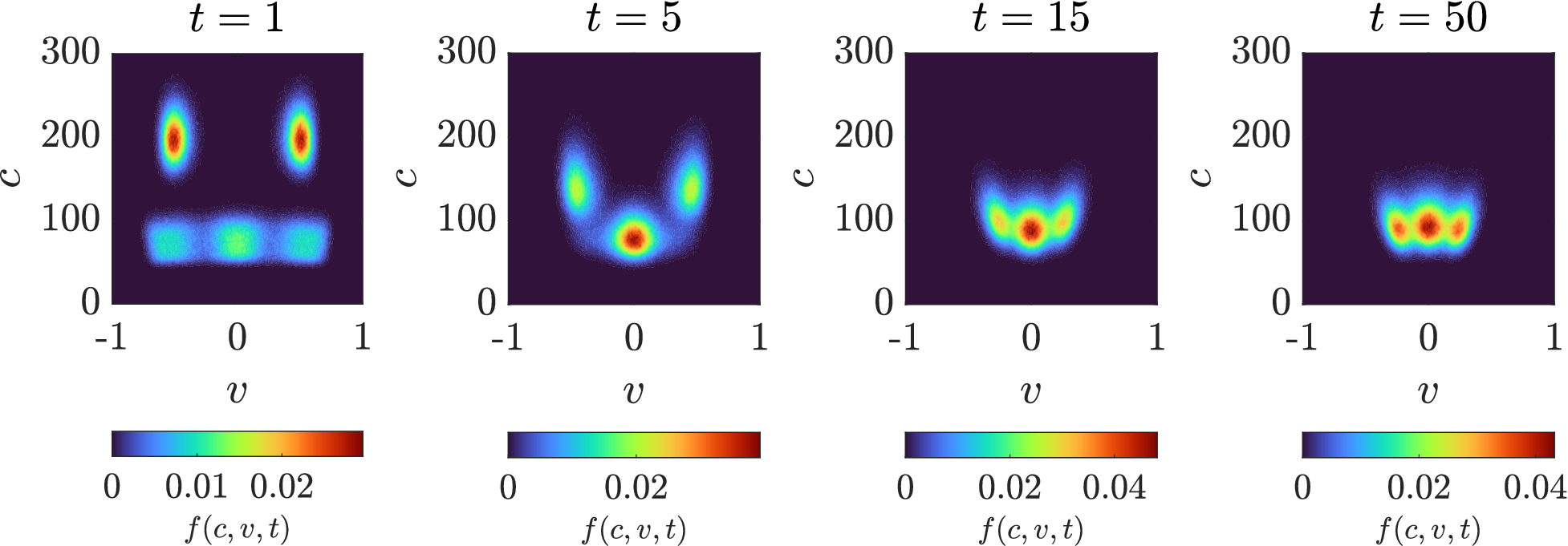}
    \caption{Control on opinions only (symmetric case).}
    \label{fig:test2_sim1}
\end{subfigure}

\par\vspace{0.5em}

\begin{subfigure}[b]{\textwidth}
    \centering
    \includegraphics[width=\textwidth]{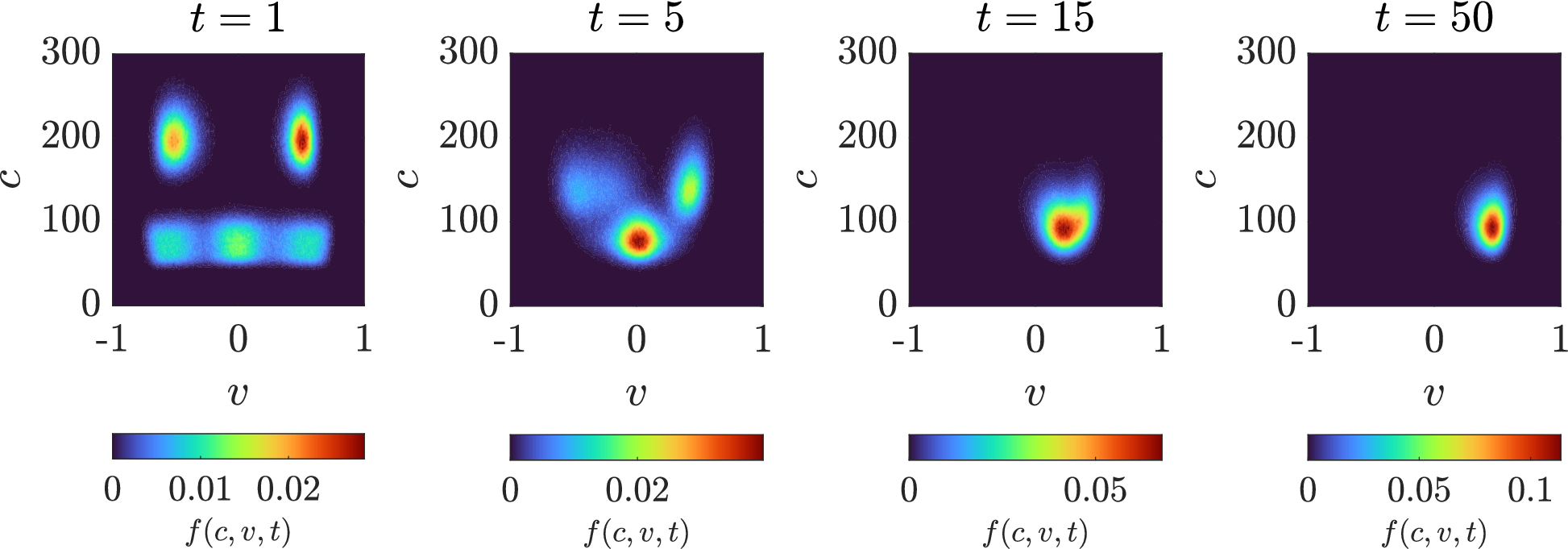}
    \caption{Control on opinions only (non-symmetric case).}
    \label{fig:test2_sim2}
\end{subfigure}

\caption{Test 2. 
Time evolution of the joint density $f(c,v,t)$ in the $(v,c)$ plane.
Each panel shows snapshots at times $t = 1, 5, 15, \revision{50}$.  
(a) uncontrolled dynamics;  
(b) symmetric opinion control;  
(c) non-symmetric opinion control.
}
\label{fig:simulations_test2}
\end{figure}

\begin{figure}[htbp]
\centering

\begin{subfigure}[b]{0.48\textwidth}
    \centering
    \includegraphics[width=\textwidth]{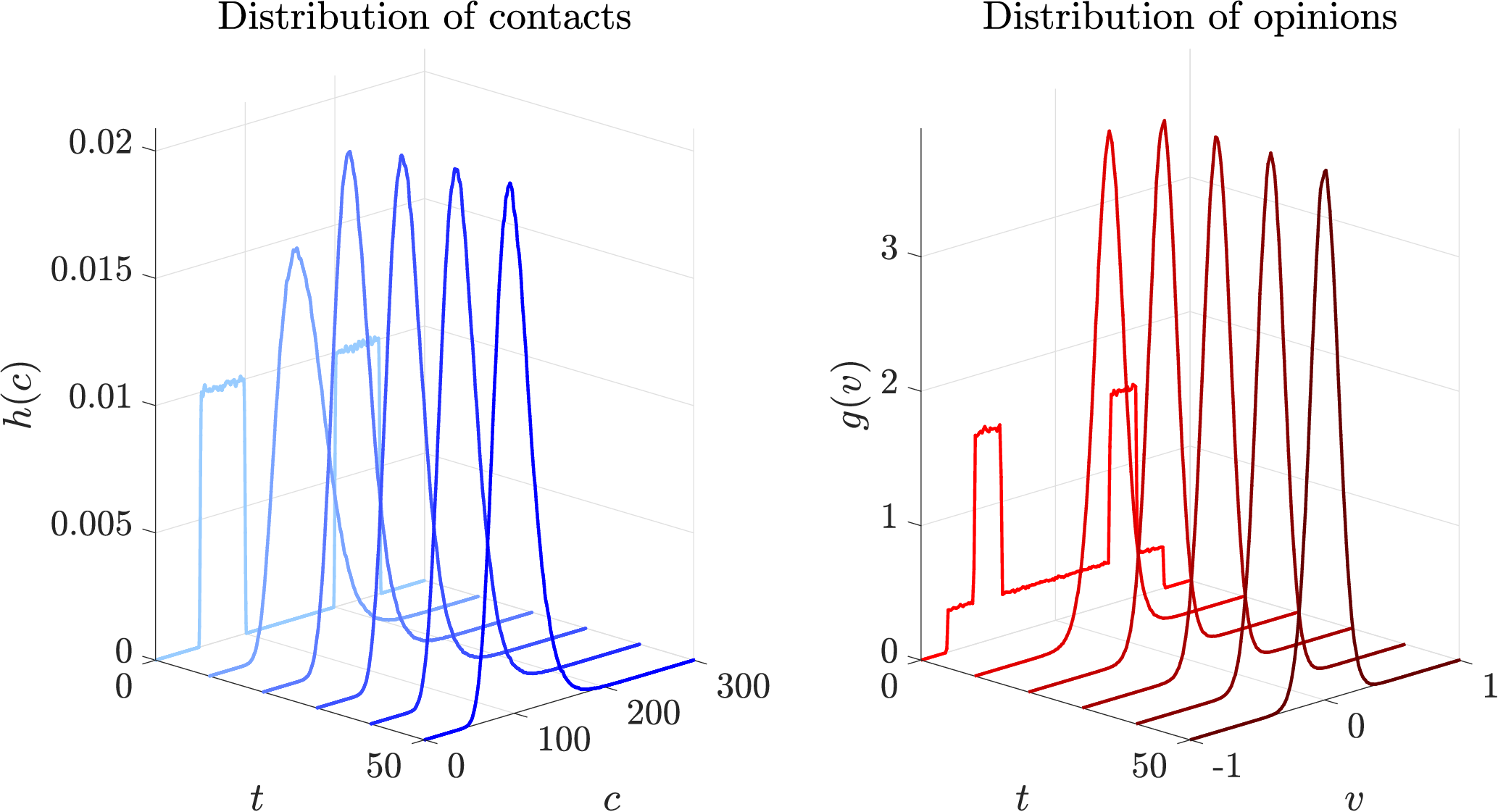}
    \caption{No control: free evolution.}
    \label{fig:dist02}
\end{subfigure}
\hfill
\begin{subfigure}[b]{0.48\textwidth}
    \centering
    \includegraphics[width=\textwidth]{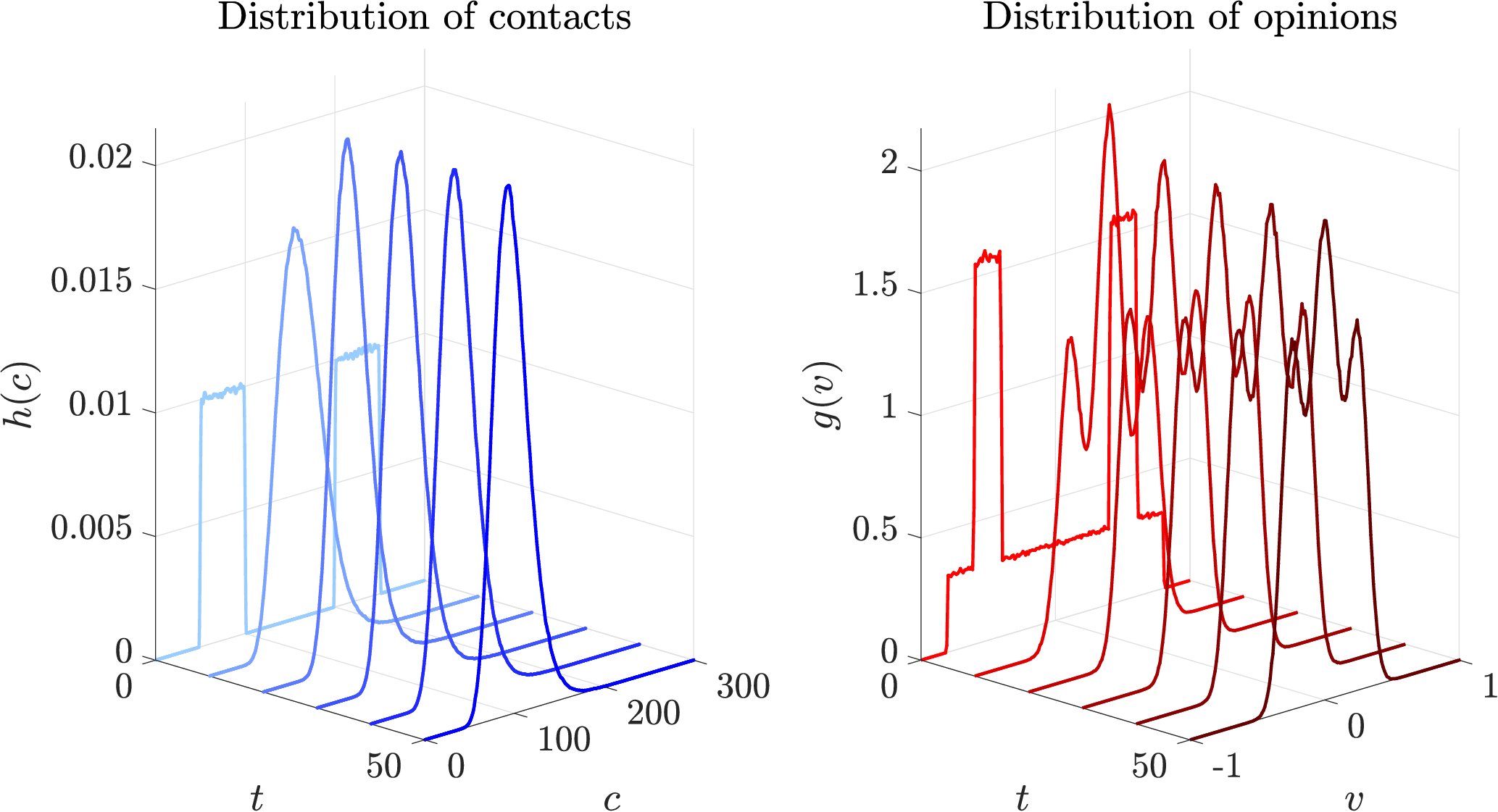}
    \caption{Symmetric opinion control.}
    \label{fig:dist12}
\end{subfigure}

\par\vspace{0.5em}

\begin{subfigure}[b]{0.48\textwidth}
    \centering
    \includegraphics[width=\textwidth]{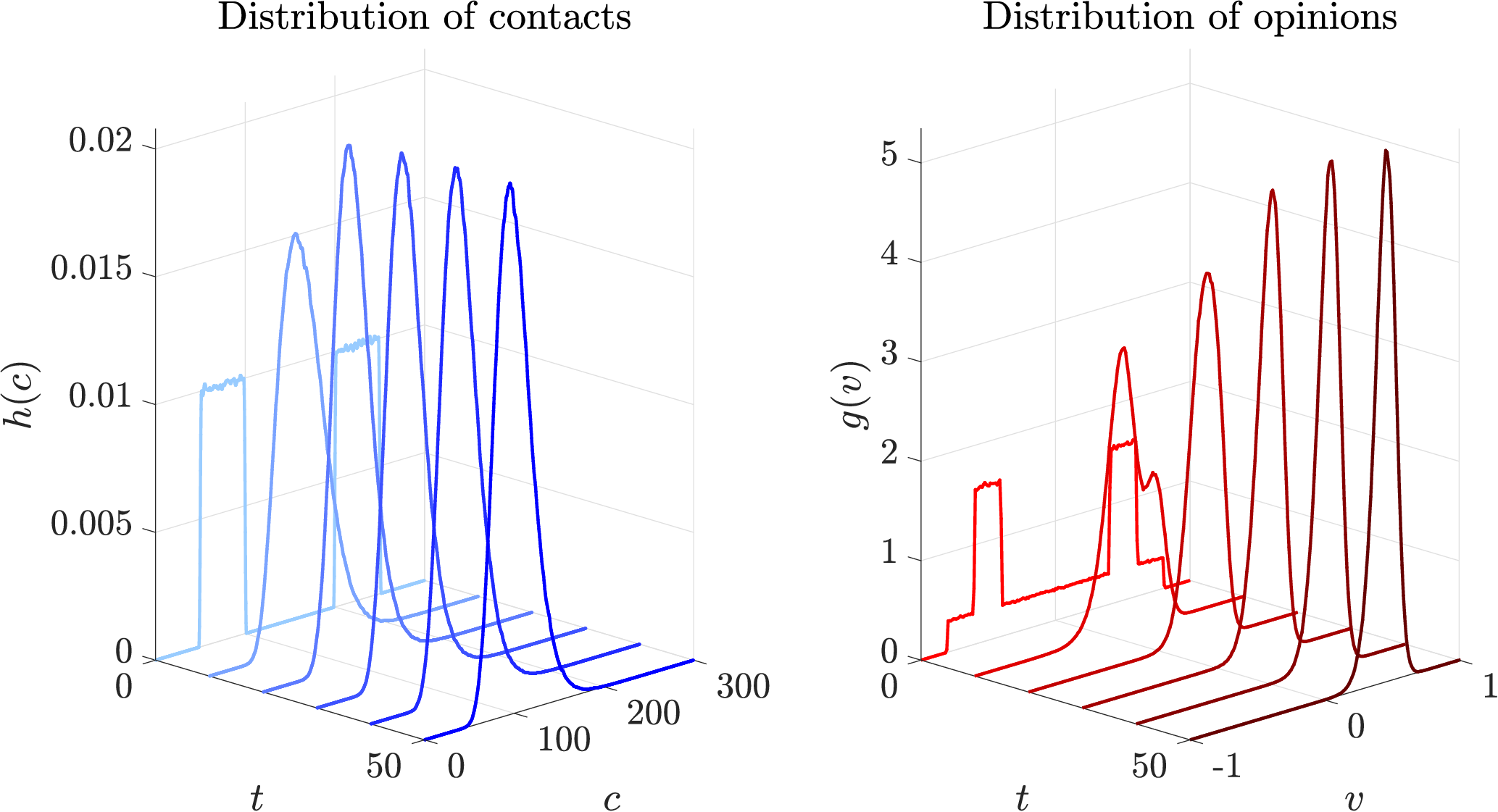}
    \caption{Non-symmetric opinion control.}
    \label{fig:dist22}
\end{subfigure}

\caption{Test 2.
Time evolution of the marginal distributions of contacts and opinions at different times. 
Each panel shows both distributions together: contacts (left) and opinions (right),
(a) uncontrolled dynamics;  
(b) symmetric opinion control;  
(c) non-symmetric opinion control.
}
\label{fig:distributions_test2}
\end{figure}

Simulation results are presented in Fig.~\ref{fig:simulations_test2}, showing snapshots of the joint density $f(c,v,t)$ at selected times, and in Fig.~\ref{fig:distributions_test2}, illustrating the evolution of the marginal distributions of contacts $h(c\revision{,t})$ and opinions $g(v\revision{,t})$.

In the uncontrolled case~(a), leaders in both groups progressively lose influence due to contact decay and are unable to maintain their initial stance. Their opinions converge toward a compromise around zero, reflecting the symmetry of the initial configuration. In scenario~(b), where both groups have the same control strength, leaders succeed in maintaining their respective target opinions. The follower population is simultaneously attracted toward both opinion poles, resulting in a persistent polarization and the emergence of a stable central cluster. Owing to symmetry, this central group remains approximately stationary and is not significantly attracted by either side. In the asymmetric scenario~(c), the stronger penalization of Group~A makes its leaders less effective at sustaining their target opinion, whereas Group~B, with a lower penalization, maintains a strong control action. Consequently, Group~B successfully drives the entire follower population toward its own target opinion, eventually dominating the competition for consensus.

In this test we see how competition between elite groups depends on the balance of resources devoted to persuasion. If neither side sustains its effort, both lose visibility and the population drifts toward a centrist consensus. Equal control efforts stabilise a polarized environment, with two enduring opinion poles and a persistent centre. When resources are asymmetric, the better-supported group dominates: unequal persuasive capacity translates into unequal influence over public opinion.

\subsection{Test 3: emergence of echo chambers under local contact control}

In this test, we analyze a setting in which two groups of leaders, characterized by markedly opposite opinions, interact with a large population of followers whose initial opinions are closer to one of the two leader groups. The aim is to study how local control of contacts performed by one or both leader groups affects the formation of echo chambers and the eventual opinion distribution of the mass.

\revision{We consider three population of interacting agents, with $25\%$ assigned to Leader Group~A ($\ell=A$), $25\%$ to Leader Group~B ($\ell=B$), and the remaining $50\%$ form the follower population ($\ell=F$).} The simulation horizon is $T=150$ with time step $\Delta t = 10^{-3}$.

\begin{figure}[htbp]
\centering
\includegraphics[width=0.4\textwidth]{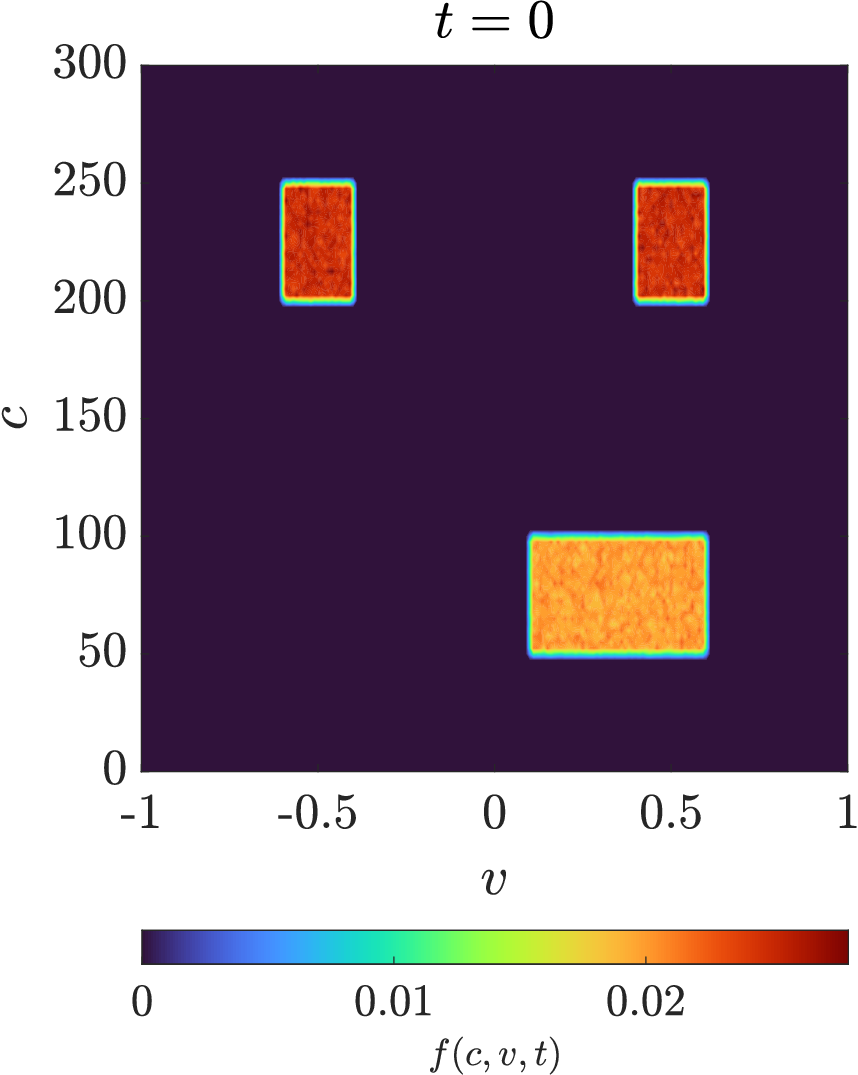}
\caption{ Test 3. Joint density $f(c,v,0)$ at the initial time.}
\label{fig:test3_t0}
\end{figure}

The dynamics of contacts and opinions evolve according to the parameters reported in Table~\ref{tab1}, here with penalization weights $\gamma_c = \gamma_v = 1$ for both populations of leaders. 

The total initial joint density $f_0(c,v)$, sampled from $N_s = 10^6$ agents, is illustrated in Fig.~\ref{fig:test3_t0}. At time $t=0$, leaders of Group~A are sampled from 
\(
(c,v)\sim \mathrm{Unif}\big([200,250]\times[-0.6,-0.4]\big),
\) and are associated with the target opinion $\tilde v_A=-0.5$, representing a highly connected community endorsing a moderately negative stance. Leaders of Group~B are initialised as 
\(
(c,v)\sim \mathrm{Unif}\big([200,250]\times[0.4,0.6]\big),
\)
with a corresponding positive target opinion $\tilde v_B = 0.5$. Followers are drawn from 
\(
(c,v)\sim \mathrm{Unif}\big([50,100]\times[0.1,0.6]\big),
\)
representing a poorly connected and more heterogeneous population.

We compare three configurations of the control dynamics. In the first one, only the opinion control $u^\star$ is active for the two leader groups. In the second configuration, both groups employ the opinion control, while Group~A also activates the contact control $\kappa^\star$. In the third configuration, both controls are active for both leader groups.

\begin{figure}[htbp]
\centering

\begin{subfigure}[b]{0.48\textwidth}
    \centering
    \includegraphics[width=\textwidth]{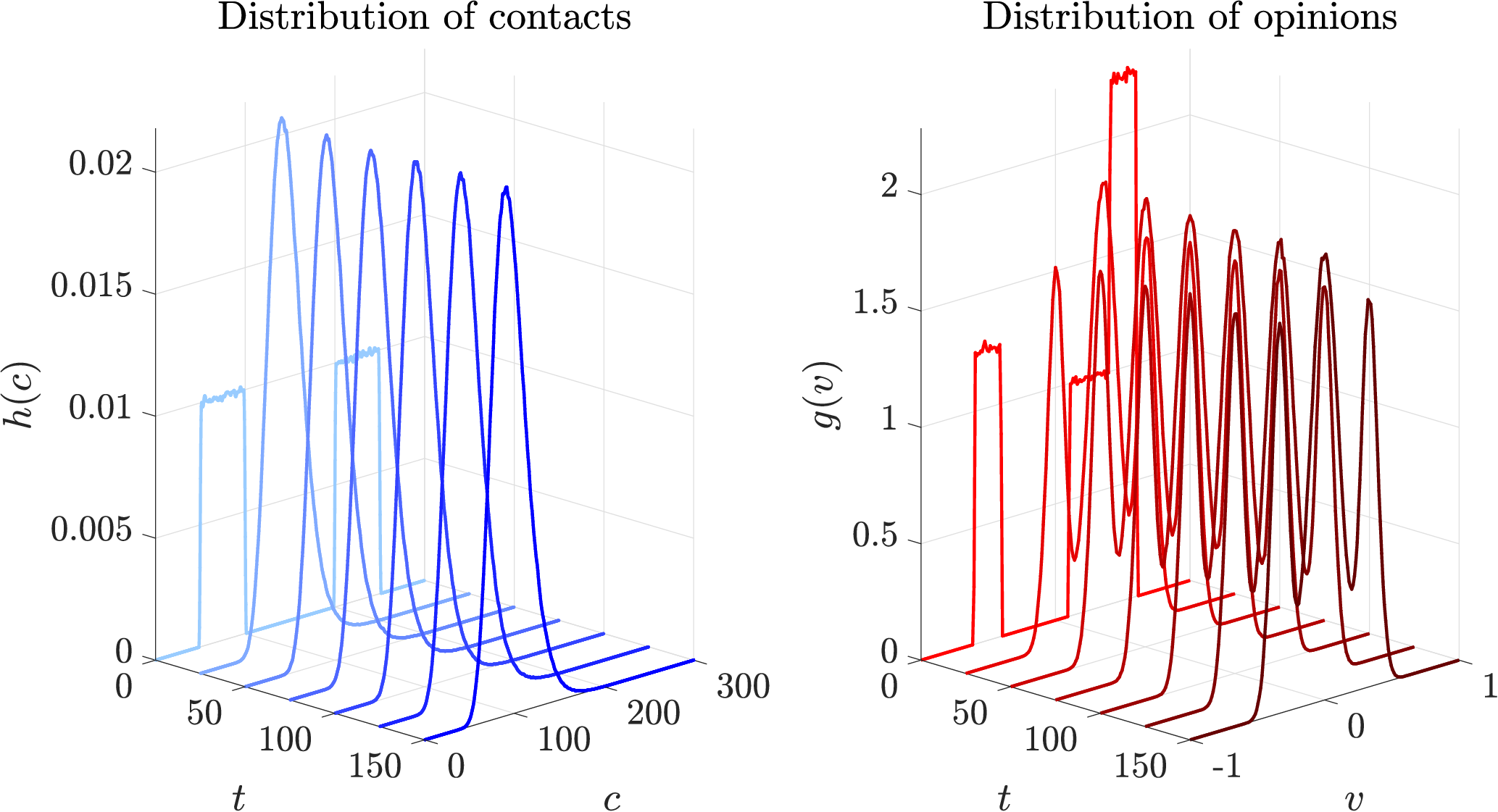}
    \caption{Control $u^\star$ on Groups A and B.}
    \label{fig:dist03}
\end{subfigure}
\hfill
\begin{subfigure}[b]{0.48\textwidth}
    \centering
    \includegraphics[width=\textwidth]{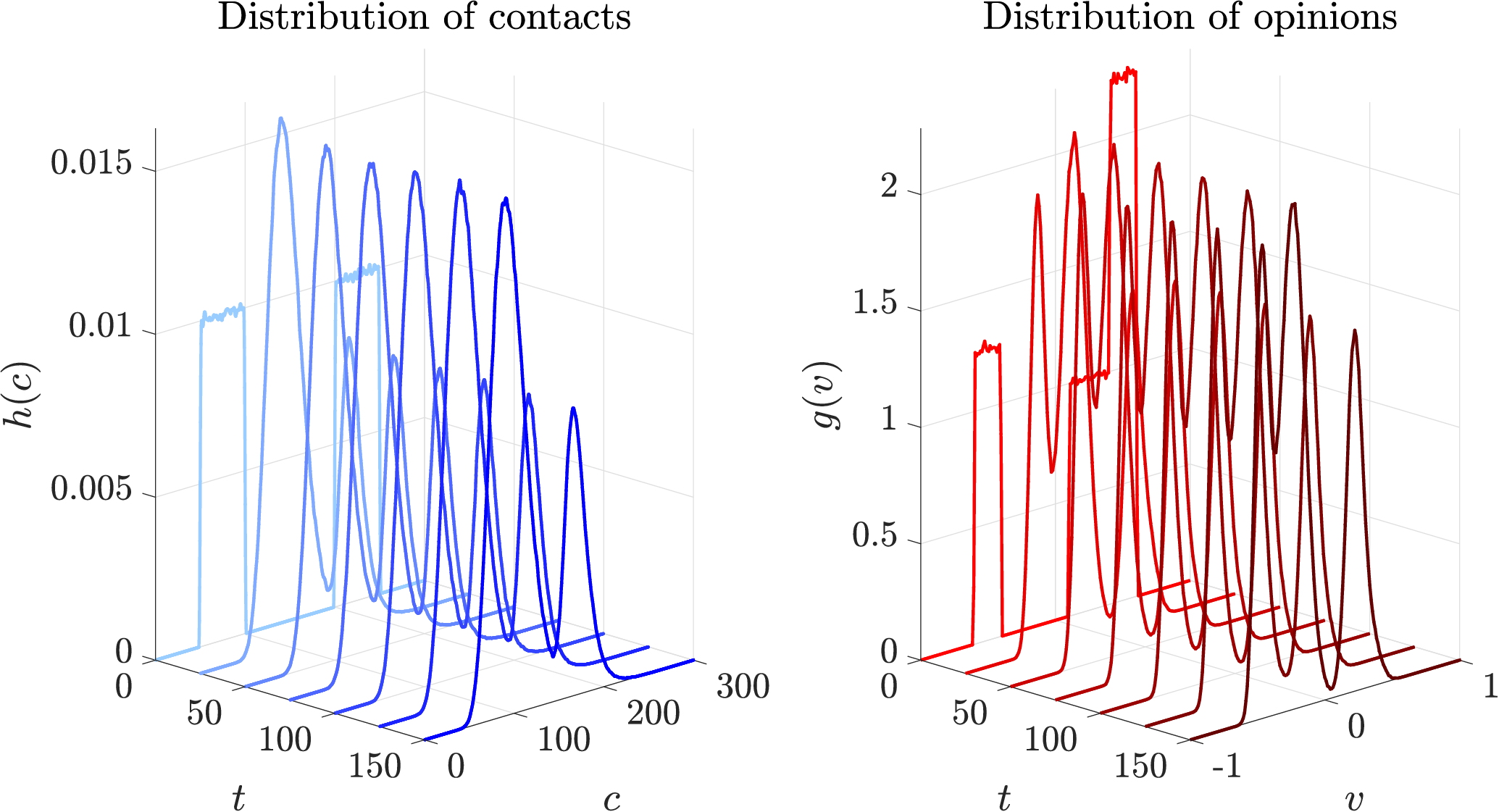}
    \caption{Control $u^\star$ on Groups A and B, control $\kappa^\star$ on Group A.}
    \label{fig:dist13}
\end{subfigure}

\par\vspace{0.5em}

\begin{subfigure}[b]{0.48\textwidth}
    \centering
    \includegraphics[width=\textwidth]{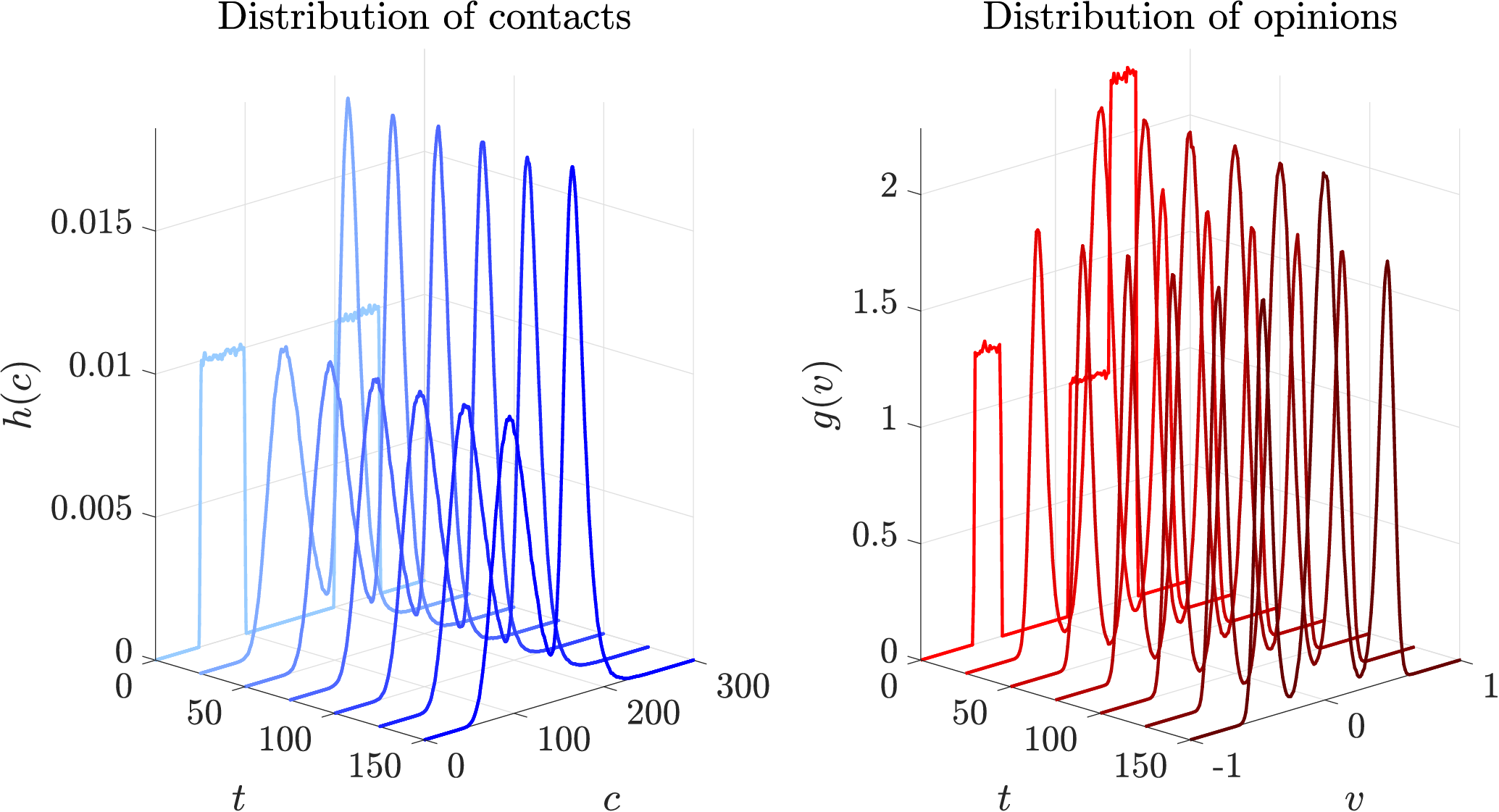}
    \caption{Controls $u^\star$ and $\kappa^\star$ on Groups A and B.}
    \label{fig:dist23}
\end{subfigure}

\caption{Test 3.
Time evolution of the marginal distributions of contacts and opinions at different times. 
Each panel shows both distributions together: contacts (left) and opinions (right),
(a) control $u^\star$ on Groups A and B;  
(b) control $u^\star$ on Groups A and B, control $\kappa^\star$ on Group A;  
(c) controls $u^\star$ and $\kappa^\star$ on Groups A and B.
}

\label{fig:distributions_test3}
\end{figure}

The results, displayed in Fig.~\ref{fig:simulations_test3} through snapshots of the total joint density $f(c,v,t)$, and in Fig.~\ref{fig:distributions_test3} via the evolution of the marginal distributions of contacts $h(c\revision{,t})$ and opinions $g(v\revision{,t})$, reveal a clear differentiation among the three scenarios. In all cases, both leader groups initially experience a decline in their number of contacts due to the endogenous mechanisms governing the contact dynamics. When only the opinion control is active, this decay leads to a strong weakening of both leader groups, and the followers tend to preserve their initial opinion distribution. When Group~A activates the contact control, however, the decline of its connectivity is mitigated and a cohesive cluster emerges around its target opinion; this cluster progressively attracts both the follower population and, indirectly, the weakened Group~B, whose loss of contacts reduces its influence. In the third scenario, where both leader groups activate the contact control, two persistent and well-connected leader clusters form around their respective target opinions, preventing their decay. In this case, a portion of the follower population converges toward an intermediate compromise opinion, remaining influenced by the simultaneous presence of the two stable leader communities: in other words, a fraction is pulled toward each pole, producing a stable bi-polarization or two echo chambers with a residual central compromise group.

This test highlights the emergence of echo chambers when leader groups preserve their visibility. Leaders who fail to maintain connectivity lose influence, while sustained visibility by one or both groups leads respectively to dominance or to the formation of two stable, polarized communities separated by a small central cluster.

\begin{figure}[htbp]
\centering

\begin{subfigure}[b]{\textwidth}
    \centering
    \includegraphics[width=\textwidth]{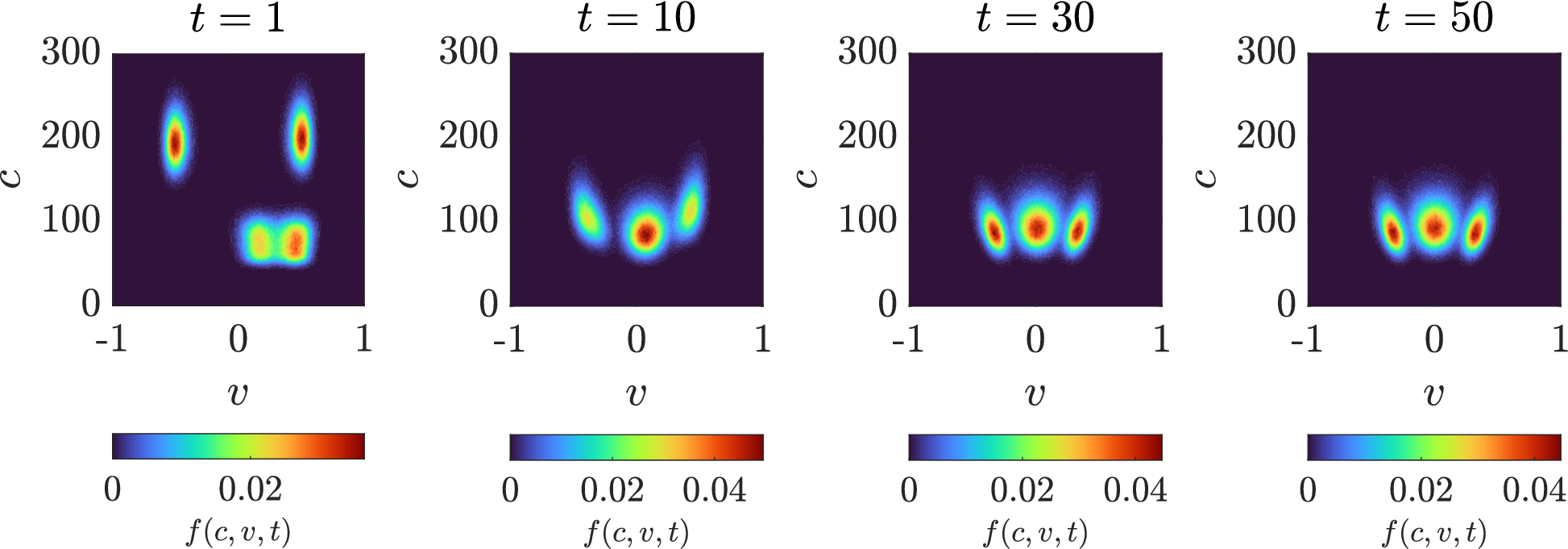}
    \includegraphics[width=\textwidth]{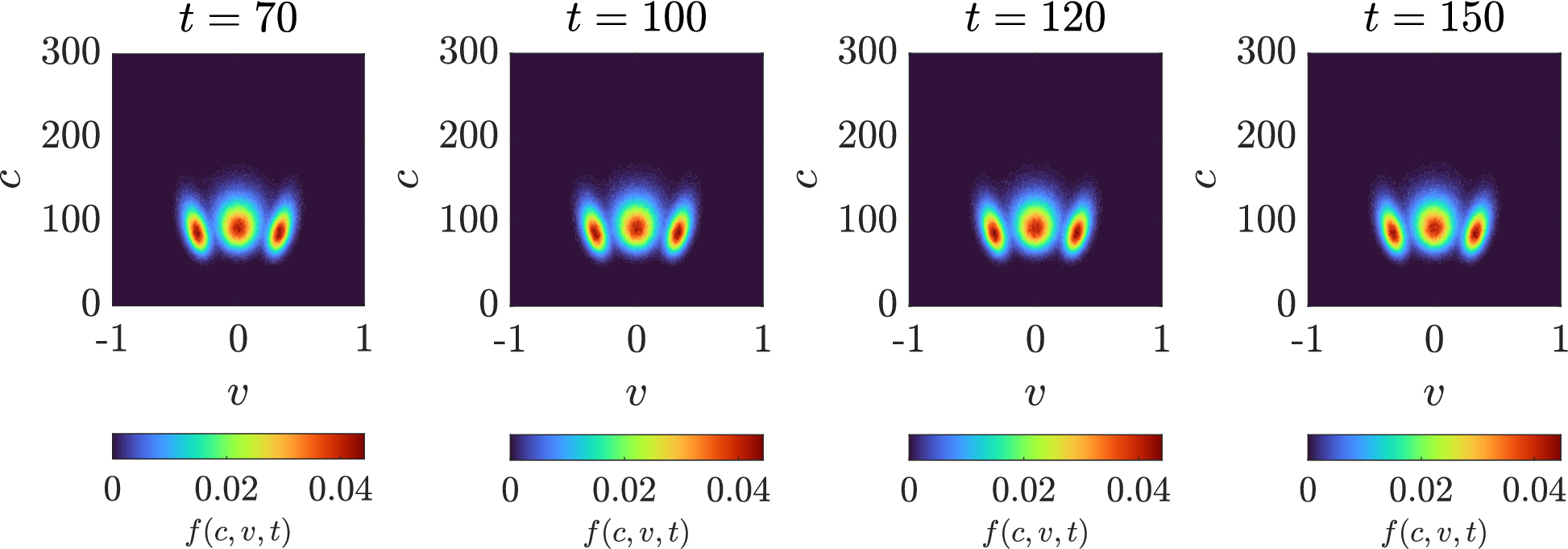}
    \caption{Control $u^\star$ on Groups A and B.}
    \label{fig:test3_sim0}
\end{subfigure}

\par\vspace{0.5em}

\begin{subfigure}[b]{\textwidth}
    \centering
    \includegraphics[width=\textwidth]{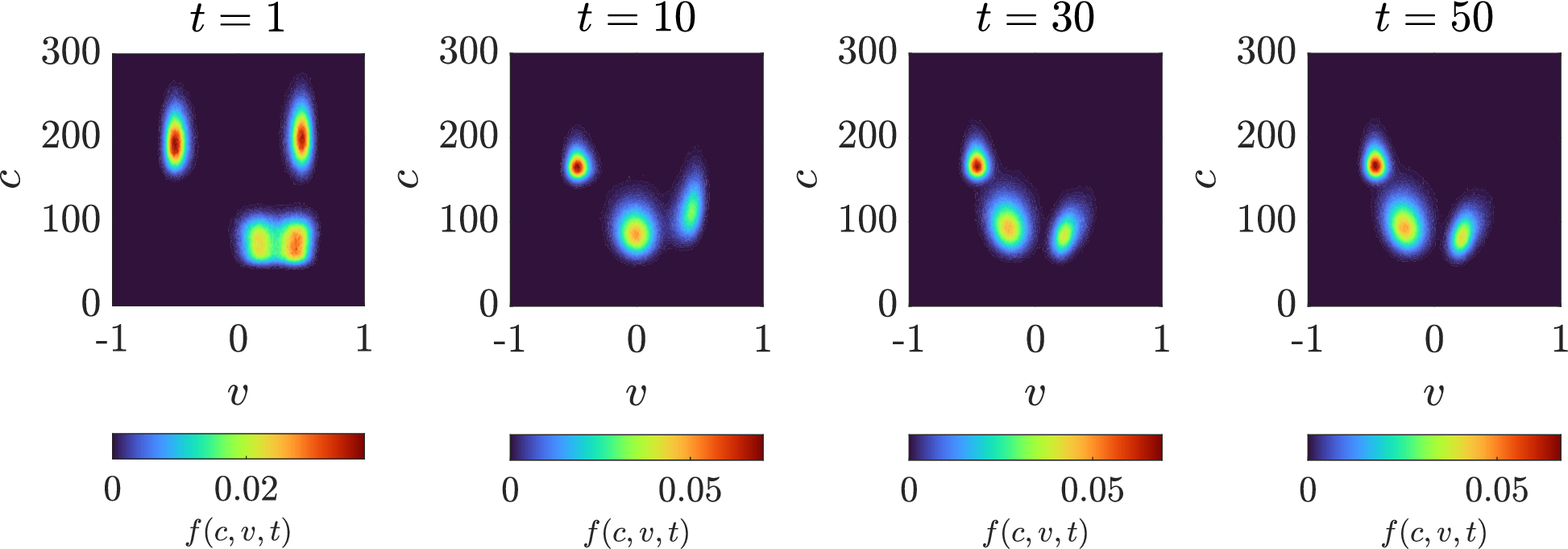}
    \includegraphics[width=\textwidth]{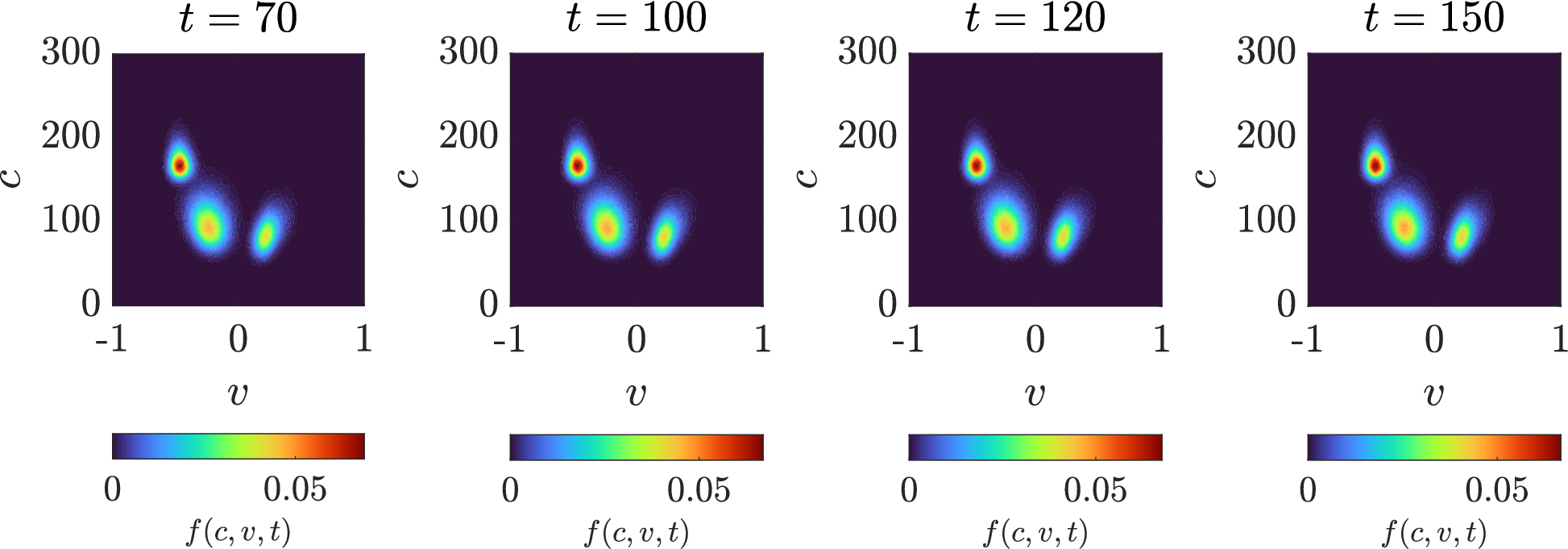}
    \caption{Control $u^\star$ on Groups A and B, control $\kappa^\star$ on Group A.}
    \label{fig:test3_sim1}
\end{subfigure}

\caption{Test 3. Time evolution of the joint density $f(c,v,t)$ in the $(v,c)$ plane.  
(a) control $u^\star$ on Groups A and B.;  
(b) control $u^\star$ on Groups A and B, control $\kappa^\star$ on Group A.}
\label{fig:simulations_test3}
\end{figure}

\clearpage

\begin{figure}[htbp]
\ContinuedFloat
\centering

\begin{subfigure}[b]{\textwidth}
    \centering
    \includegraphics[width=\textwidth]{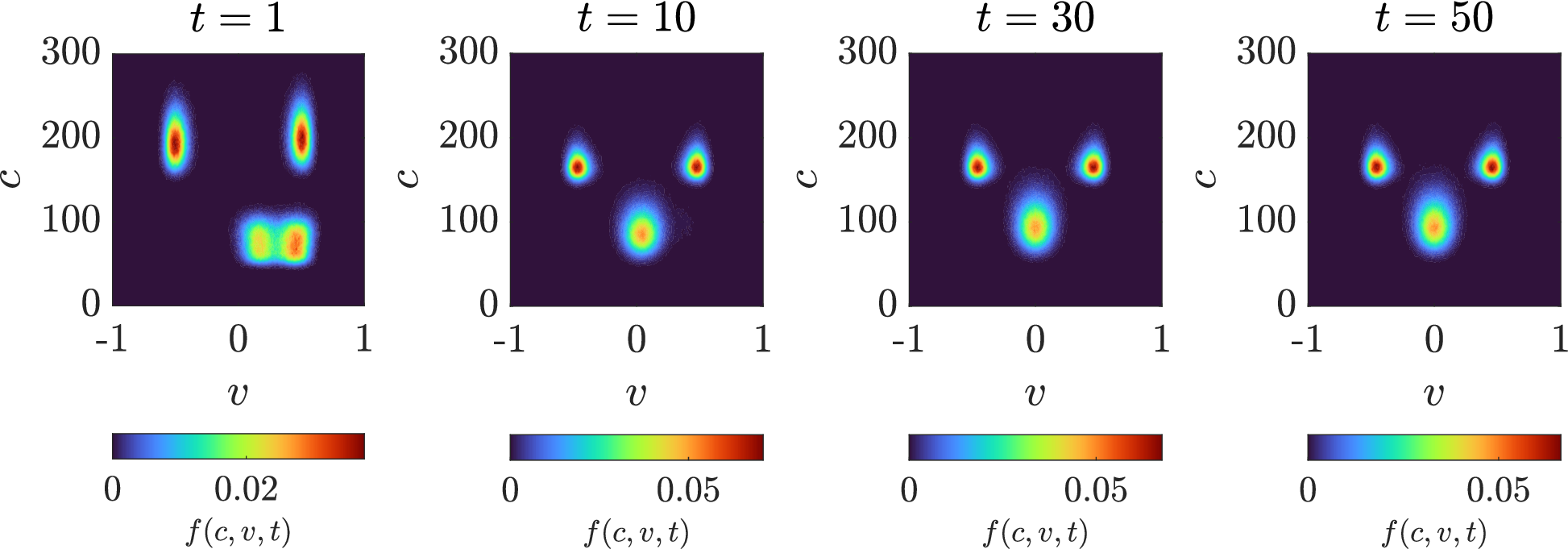}
    \includegraphics[width=\textwidth]{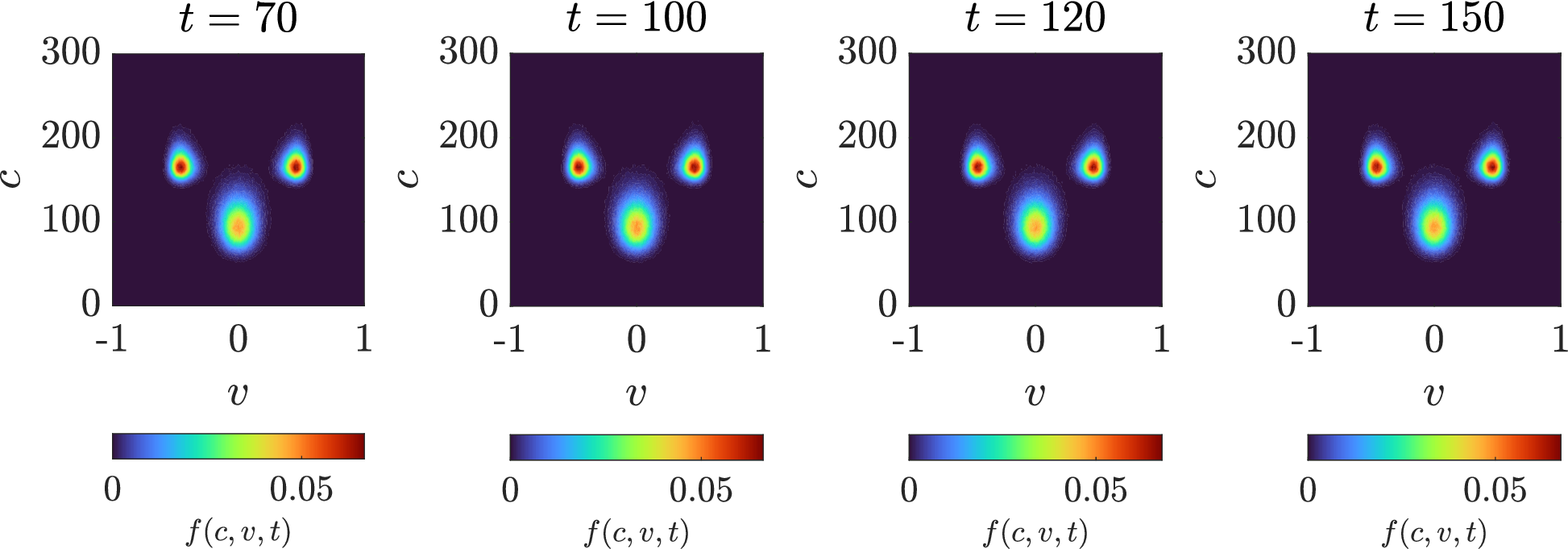}
    \caption{Controls $u^\star$ and $\kappa^\star$ on Groups A and B.}
    \label{fig:test3_sim2}
\end{subfigure}

\caption{(c) controls $u^\star$ and $\kappa^\star$ on Groups A and B.}
\end{figure}

\section{Conclusions}
In this paper, we have introduced a kinetic modeling framework for opinion dynamics in popularity-adaptive social networks, in which the evolution of opinions is explicitly coupled with the formation and adaptation of social connections. A central feature of the proposed approach is the presence of a feedback mechanism where agents’ opinions influence the dynamics of their social contacts, describing the interplay between influence and connectivity.
Furthermore, the model include control strategies acting simultaneously on popularity dynamics and opinion exchanges. We believe that such modelling framework allows to reproduce salient features of online social systems, including visibility effects, leader–follower interactions, and opinion-driven reorganization of network connections. 
Possible extensions of the present work include data-driven calibration and validation of the model using empirical social network data, and the formulation of strategic, game-theoretic settings involving multiple competing controllers with conflicting influence objectives.

\section*{Acknowledgments}
This work has been written within the activities of GNCS groups of INdAM
(Italian National Institute of High Mathematics). GA and MP have been partially supported by
MUR-PRIN Project 2022 No. 2022N9BM3N “Efficient numerical schemes and optimal
control methods for time-dependent partial differential equations” financed by the European Union Next Generation EU. GA, EC, and MP have been partially supported by MUR–PRIN 2022 through the PNRR Project No. P2022JC95T “Data-driven
discovery and control of multi-scale interacting artificial agent systems” financed by the European Union Next Generation EU. EC acknowledges the support by Fondo Italiano per la Scienza (FIS2023-01334) advanced grant "ADvanced numerical Approaches for MUltiscale Systems with uncertainties" - ADAMUS.
\bibliographystyle{plain}
\bibliography{biblio_opinion}

\end{document}